\documentclass{amsart}
\usepackage[all]{xy}
\usepackage{amsfonts,amssymb,amsmath,amscd}
\usepackage{amsthm}

\newcounter{zlist}
\newenvironment{zlist}{\begin{list}{{\rm(\arabic{zlist})}}{
\usecounter{zlist}\leftmargin2.5em\labelwidth2em\labelsep0.5em
\topsep0.6ex\itemsep0.3ex plus0.2ex minus0.3ex
\parsep0.3ex plus0.2ex minus0.1ex}}{\end{list}}

\newcounter{blist}
\newenvironment{blist}{\begin{list}{{\rm(\alph{blist})}}{
\usecounter{blist}\leftmargin2.5em\labelwidth2em\labelsep0.5em
\topsep0.6ex \itemsep0.3ex plus0.2ex minus0.3ex
\parsep0.3ex plus0.2ex minus0.1ex}}{\end{list}}

\newcounter{rlist}
\newenvironment{rlist}{\begin{list}{{\rm(\roman{rlist})}}{
\usecounter{rlist}\leftmargin2.5em\labelwidth2em\labelsep0.5em
\topsep0.6ex\itemsep0.3ex plus0.2ex minus0.3ex
\parsep0.3ex plus0.2ex minus0.1ex}}{\end{list}}

\swapnumbers
\newtheorem{theorem}{Theorem}[section]
\newtheorem{lemma}[theorem]{Lemma}
\newtheorem{thm}[theorem]{}
\newtheorem{proposition}[theorem]{Proposition}
\newtheorem{definition}[theorem]{Definition}

\numberwithin{equation}{section}

\newcommand{\A}{{\mathbb{A}}}
\newcommand{\B}{{\mathbb{B}}}
\newcommand{\M}{{\mathbb{M}}}
\newcommand{\II}{{\mathbb{I}}}
\newcommand{\V}{{\mathbb{V}}}
\newcommand{\C}{{\mathbb{C}}}
\newcommand{\X}{{\mathbb{X}}}

\newcommand{\cV}{{\mathcal{V}}}

\newcommand{\uA}{\mathbf{A}}
\newcommand{\oC}{\mathbf{C}}
\newcommand{\uC}{\mathbf{C}}

\newcommand{\ot}{\otimes}
\newcommand{\R}{{\rm Rat}}

\newcommand{\Hom}{{\rm Hom}}
\newcommand{\Nat}{{\rm Nat}}
\newcommand{\Set}{{\rm Set}}
\newcommand{\ev}{{\rm ev}}
\newcommand{\di}{\diamond}

\newcommand{\id}{I}
\newcommand{\lra}{\longrightarrow}

\newcommand{\Ra}{\Rightarrow}
\newcommand{\bC}{\mathbf{C}}
\newcommand{\bG}{\mathbf{G}}

\newcommand{\bT}{\mathbf{T}}
\newcommand{\rG}{G}
\newcommand{\rH}{H}
\newcommand{\rT}{T}
\newcommand{\ve}{\varepsilon}
\newcommand{\pp}{\mathcal{P}}

\begin{document}

\title{On Rational Pairings of Functors}
\author{Bachuki Mesablishvili and Robert Wisbauer}

\begin{abstract}
In the theory of coalgebras $C$ over a ring $R$, the rational
functor relates the category of modules over the algebra $C^*$
(with convolution product) with the category of comodules over $C$.
It is based on the pairing of the algebra $C^*$ with the coalgebra $C$
  provided by the evaluation map $\ev:C^*\ot_R C\to R$.

We generalise this situation by defining a {\em pairing}
between endofunctors $T$ and $G$ on any category $\A$ as
a map, natural in $a,b\in \A$,
  $$\beta_{a,b}:\A(a, G(b)) \to \A(T(a),b),$$
and we call it {\em rational} if these all are injective.
In case $\bT=(T,m_T,e_T)$ is a monad and $\bG=(G,\delta_G,\ve_G)$
is a comonad on $\A$,
additional compatibility conditions are imposed on a pairing between $\bT$
and $\bG$. If such a pairing is given and is rational, and
 $\bT$ has a right adjoint monad $\bT^\di$,
we construct a {\em rational functor} as the functor-part of an
idempotent comonad on the $\bT$-modules $\A_{\rT}$ which
generalises the crucial properties of the rational functor for
coalgebras. As a special case we consider pairings on monoidal
categories.
\end{abstract}

\maketitle

\tableofcontents

\section{Introduction}

The pairing of a $k$-vector space $V$ with its dual space $V^*=\Hom(V,k)$
provided by the evaluation map $V^*\ot V\to k$ can be extended
from base fields $k$ to arbitrary base rings $A$. Then it can be applied to the study of $A$-corings $C$ to obtain a faithful functor from the
category of $C$-comodules to the category of $C^*$-modules.
The purpose of this paper it to extend these results to (endo)functors on
arbitrary categories.
We begin by recalling some facts from module theory.

\begin{thm}\label{Example}{\bf Pairing of modules.} \em
Let $C$ be a bimodule over a ring $A$ and
 $C^*=\Hom_A(C,A)$ the right dual.
Then $C\ot_A-$ and $C^*\ot_A-$ are endofunctors
on the category $_A\M$ of left $A$-modules
and the evaluation
$$\ev: C^*\ot_A C \to A, \quad f\ot c\mapsto f(c),$$
induces a pairing between these functors.
For left $A$-modules $X,Y$, the map
$$\alpha_Y :C\ot_A Y\to {_A\Hom(C^*,Y)}, \quad c\ot y \mapsto [f\mapsto f(c)y],$$
induces the map
  $$\begin{array}{rcl}
\beta_{X,Y}:
  {_A\Hom (X, C\ot_A Y)}  & \lra & {_A\Hom( X, {_A\Hom(C^*,Y)})}, \\[+1mm]
  X\stackrel{f}\to C\ot_A Y & \longmapsto &  X\stackrel{f}\to C\ot_A Y
   \stackrel{\alpha_Y}\lra {_A\Hom(C^*,Y)}.
\end{array} $$
Clearly $\beta_{X,Y}$ is injective
for all left $A$-modules $X,Y$ if and only if $\alpha_Y$
is a monomorphism (injective) for any left $A$-module $Y$, that is,
$C_A$ is locally projective (see \cite{Abu}, \cite[42.10]{BW}).
 \end{thm}

Now consider the situation above with some additional structure.

\begin{thm}\label{pair-cor} {\bf Pairings for corings.} \em
Let ${\bC}= (C,\Delta,\ve)$ be a coring over the ring $A$, that is,
$C$ is an $A$-bimodule with bimodule morphisms
coproduct $\Delta:C\to C\ot_A C$ and counit $\ve:C\to A$.
Then the right dual $C^*=\Hom_A(C,A)$
has a ring structure by the convolution product for $f,g\in {C}^*$,
$f*g= f \circ(g \ot_A \id_C) \circ\Delta$
(convention opposite to \cite[17.8]{BW}) with
unit $\ve$, and we have a pairing between the comonad $C\ot_A-$
and the monad $C^*\ot_A-$ on $_A\M$. In this case,
${\Hom_A(C^*,-)}$ is a comonad on $_A\M$ and $\alpha_Y$ considered in
\ref{Example} induces a  comonad morphism
$\alpha:C\ot-\to \Hom_A(C^*,-)$. We have the commutative diagrams
\begin{equation}\label{coring}
\xymatrix{C^*\ot_A C^*\ot_A C\ar[r]^{\id\ot\id\ot \Delta\quad\;}
 \ar[d]_{*\ot \id} & C^*\ot_A C^*\ot_A C\ot_AC \ar[r]^{\qquad\;\;\id\ot \ev\ot\id}
 &  C^*\ot_A C \ar[d]^\ev & \ar[l]_{\quad\ve\ot \id } C \ar[dl]^\ve\\
  C^*\ot_A C \ar[rr]^\ev && A & .}
\end{equation}

The Eilenberg-Moore category  $\M^{\Hom(C^*,-)}$ of ${\Hom(C^*,-)}$-comodules is equivalent to the category ${_{C^*}\M}$ of left
$C^*$-modules (e.g. \cite[Section 3]{BBW}) and thus $\alpha$ induces a functor
 $${^C\M}\to \M^{\Hom(C^*,-)} \simeq {_{C^*}\M} $$
which is fully faithful if and only if
the pairing $(C^*,C,\ev)$ is
rational, that is $\alpha_Y$ is monomorph for all $Y\in {_A\M}$
(see \cite[19.2 and 19.3]{BW}).
Moreover, $\alpha$ is an isomorphism if and only if the categories
$^C\M$ and $_{C^*}\M$ are equivalent and this is tantamount to
$C_A$ being finitely generated and projective.
\end{thm}

In Section 2 we recall the notions and some basic facts on natural transformations 
between endofuctors needed for our investigations.

Weakening the conditions for an adjoint pair of functors, a {\em pairing}
 of two functors $T:\A\to \B$ and $G:\B\to \A$ is defined as a map
$\beta_{a,b}:\A(a, G(b)) \to \A(T(a),b)$, natural in $a\in \A$, $b\in \B$
of two functors between arbitrary categories is defined
in Section 3
(see \ref{pair-func}) and it is called {\em rational} if
all the $\beta_{a,b}$ are injective maps.
For pairing of monads $\bT$ with comonads $\bG$ on a category $\A$,
additional conditions
are imposed on the defining natural transformations (see \ref{pair-mon}).
These imply the existence of a functor $\Phi^\pp: \A^{\rG} \to \A_T$
from the $\bG$-comodules to the $\bT$-modules (see \ref{func-ind}),
 which is full and faithful provided the pairing is rational (see \ref{rat-pair}).
Of special interest is the situation that the monad $\bT$ has a right adjoint
$\bT^\di$ and the last part of Section 3 is dealing with this case.

  Referring to these results, a {\em rational functor} $\R^\pp:\A_T\to \A_T$
is associated with any rational pairing in Section 4. This leads to the definition of {\em rational} $T$-modules and under some additional conditions they form a
coreflective subcatgeory of $\A_T$ (see \ref{T.2.15}).

The application of the general notions of pairings to monoidal categories
is outlined in Section 5. The resulting formalism is very close to the
module case considered in \ref{pair-cor}.

In Section 6, we apply our results to entwining structures $(\uA,\uC,\lambda)$
 on monoidal categories $\cV=(\V, \ot, \II)$. 
The objects $A$ and $C$ incude a functor $\V(-\ot \!C, A):\V^{op}\to \rm{Set}$
and if this is representable we call the entwining
\emph{representable}, that is if $\V(-\ot C, A)\simeq \V(-,E)$ for some 
 object $E\in \V$. This $E$ allows for an algebra 
structure, an algebra morphism $A\to E$ (see Proposition \ref{attached-alg})
and a functor from the category ${^C_A\V}(\lambda)$  of entwined modules to the category ${_E\V}$ of left $E$-modules. In case the tensor functors have right
 adjoints a pairing on $\V$ is related to the entwining (see \ref{p-lambda}) and
 its properties are studed. 
Several results known for the rational functors for ordinary entwined modules (see, for example, \cite{Abu}, \cite{EG} and \cite{EGL} )
can be obtained as corollaries from the main result of this section 
(see Theorem \ref{th-entw}).

\section{Preliminaries}

In this section we recall some notation and basic facts from category theory.
Throughout $\A$ and $\B$ will denote any categories. By $\id_a$, $\id_\A$ or
just by $\id$ we denote the identity morphism of an object $a\in \A$,
respectively the identity functor of a category $\A$.

\begin{thm}\label{mon-comon}{\bf Monads and comonads.} \em
For a monad $\bT=(T, m_{T}, e_{T})$ on
$\A$, we write
\begin{rlist}
\item[] $\A_{\rT}$ for the Eilenberg-Moore category of $\bT
$-modules;
\item[] $U_{\rT} : \A_{\rT} \to \A ,\,\, (a, h_a) \to a,$ for the
underlying (forgetful) functor;
\item[] $\phi_{\rT} : \A \to \A_{\rT},\,\, a \to (T(a),(m_{T})_a),$ for the free
$\bT$-module functor, and
\item[] $\eta_{\rT}, \varepsilon_{\rT}: \phi_{\rT}
\dashv U_{\rT} : \A_{\rT} \to \A$ for the
forgetful-free adjunction.
\end{rlist}

Dually, if ${\bG} =(G, \delta_{G},
\varepsilon_{G})$ is a comonad on $\A$,  we write
\begin{rlist}
\item[] $\A^{\rG}$ for the category of the Eilenberg-Moore category of
${\bG}$-comodules;
\item[] $U^{\rG}: \A^{\rG} \to \A, \,\, (a, \theta_a) \to a$, for the
forgetful functor;
\item[] $\phi^{\rG} : \A \to \A^{\rG},\,\, a \to (G(a), (\delta_{G})_a)$, for
the cofree ${\bG}$-comodule functor, and
\item[] $\eta^{\rG} , \varepsilon^{\rG} : U^{\rG} \dashv \phi^{\rG} : \A \to \A^{\rG}$ for the
forgetful-cofree adjunction.
\end{rlist}
\end{thm}

\begin{thm}\label{idempotent}{\bf Idempotent comonads.} \em A
comonad $\textbf{H}=(H,\varepsilon,\delta)$ is said to be
\emph{idempotent} if one of the following equivalent conditions
is satisfied (see, for example, \cite{CW}):
\begin{blist}
    \item  the forgetful functor $U^H : \A^H \to \A$ is full
    and faithful;
    \item the unit $\eta^H:\id\to\phi^H U^H$ of the adjunction
      $U^H \dashv \phi^H$ is an isomorphism;
    \item  $\delta : H \to HH$ is an isomorphism;
    \item  for any $(a,\vartheta_a) \in \A^H$, the morphism
    $\vartheta_a : a \to H(a)$ is an isomorphism and $(\vartheta_a)^{-1}=\ve_a$;
    \item $H\varepsilon$ (or $\varepsilon H$) is an isomorphism.
\end{blist}

When $\textbf{H}$ is an idempotent comonad, then an object $a \in
\A$ is the carrier of an $H$-comodule if and only if there exists
an isomorphism $H(a) \simeq a$, or, equivalently, if and only if
the morphism $\varepsilon_a : H(a) \to a$ is an isomorphism. In
this case, the pair $(a, (\varepsilon_a)^{-1})$ is an
$H$-comodule. In fact, every $H$-comodule is of this form. Thus
the category $\A^\rH$ is isomorphic to the full subcategory of
$\A$ generated by those objects for which there exists an
isomorphism $H(a) \simeq a$.

In particular, any comonad $\textbf{H}=(H,\varepsilon,\delta)$
 with $\varepsilon$ a componentwise monomorphism is
idempotent (e.g. \cite{CW}). In this case, there is at most one
morphism from any comonad $\textbf{H}'$ to $\textbf{H}$. When
$\textbf{H}'$ is also an idempotent comonad, then a natural
transformation $\tau: \textbf{H}' \to \textbf{H}$ is a morphism of
comonads if and only if $\varepsilon \cdot \tau=\varepsilon',$
where $\varepsilon'$ is the counit of the comonad $\textbf{H}'$.
\end{thm}

We will need the following result whose proof is an easy diagram
chase:

\begin{lemma} \label{2.2.1} Suppose that in the commutative diagram
$$\xymatrix{a \ar[d]_{h_1}\ar[r]^{k}&b \ar[d]_{h_2}
\ar@<-0.5ex>[r]_-{g}
\ar@<0.5ex>[r]^-{f} & c \ar[d]^{h_3}\\
a' \ar[r]_{k'}& b'\ar@<-0.5ex>[r]_-{g'} \ar@<0.5ex>[r]^-{f'} &
c',}$$
the bottom row is an equaliser and the morphism $h_3$ is a
monomorphism. Then the left square in the diagram is a pullback if
and only if the top row is an equaliser diagram.
\end{lemma}

\begin{proposition}\label{p.2.2}
Let $t: G \to R $ be a natural transformation between endofunctors of $\A$
with componentwise monomorphisms
and assume that $R$ preserves equalisers. Then the following are equivalent:

\begin{blist}
\item the functor $G$ preserves equalisers;
\item  for any regular monomorphism $i: a_0 \to a$ in $\A$, the following square
  is a pullback:
$$
\xymatrix{G(a_0) \ar[r]^{G(i)} \ar[d]_{t_{a_0}} & G(a)\ar[d]^{t_a}\\
R(a_0) \ar[r]_{R(i)} & R(a).}$$
\end{blist}
 When $\A$ admits and $G$ preserves pushouts, \emph{(a)} and
\emph{(b)} are equivalent to:

\begin{blist}
\setcounter{blist}{2} \item $G$ preserves regular monomorphisms,
i.e. $G$ takes a regular monomorphism into a regular monomorphism.
\end{blist}
\end{proposition}
\begin{proof}
(b)$\Ra$(a) Let
$$
\xymatrix{a \ar[r]^{k}&b \ar@<-0.5ex>[r]_-{g} \ar@<0.5ex>[r]^-{f}
& c }$$ be an equaliser diagram in $\A$ and consider the
commutative diagram
$$
\xymatrix{G(a) \ar[d]_{t_a}\ar[r]^{G(k)}& G(b) \ar[d]_{t_b}
\ar@<-0.5ex>[rr]_-{G(g)}
\ar@<0.5ex>[rr]^-{G(f)} && G(c) \ar[d]^{t_c}\\
R(a) \ar[r]_{R(k)}& R(b) \ar@<-0.5ex>[rr]_-{R(g)}
\ar@<0.5ex>[rr]^-{R(f)} && R(c)\,.}$$ Since $R$ preserves
equalisers, the bottom row of this diagram is an equaliser. By
(b), the left square in the diagram is a pullback. $t_c$ being a
monomorphism, it follows from Lemma \ref{2.2.1} that the top row
of the diagram is an equaliser. Thus $G$  preserves equalisers.
\smallskip

(a)$\Ra$(b) Reconsider the diagram in the above proof and apply Lemma \ref{2.2.1}.
\smallskip

Suppose now that $\A$ admits and $G$ preserves pushouts.
The implication (a)$\Ra$(c) always holds and so it remains to show

 (c)$\Ra$(b) Consider an arbitrary regular
monomorphism $i: a_0 \to a$ in $\A$. Since $\A$ admits pushouts
and $i$ is a regular monomorphism in $\A$, the diagram
$$
\xymatrix{a_0 \ar[r]^{i}& a \ar@<-0.5ex>[r]_-{i_2}
\ar@<0.5ex>[r]^-{i_1} & a \!\sqcup_{a_0} \!a\, ,}$$ where $i_1$
and $i_2$ are the canonical injections into the pushout, is an
equaliser diagram (e.g. \cite[Proposition 11.22]{AR}).
Consider now the commutative diagram
$$
\xymatrix{G(a_0) \ar[d]_{t_{a_0}}\ar[rr]^{G(i)}&& G(a)
\ar[d]_{t_a} \ar@<-0.5ex>[rr]_-{G(i_2)}
\ar@<0.5ex>[rr]^-{G(i_1)} && G(a\! \sqcup_{a_0} \!a) \ar[d]^{t_{a \sqcup_{a_0} \!a}}\\
R(a_0) \ar[rr]^{R(i)}& & R(a) \ar@<-0.5ex>[rr]_-{R(i_2)}
\ar@<0.5ex>[rr]^-{R(i_1)} && R(a \!\sqcup_{a_0} \!a) \,,}$$ in
which the bottom row is an equaliser diagram since $R$ preserves
equalisers. Since $G$ takes regular monomorphisms into regular
monomorphisms and $G$ preserves pushouts, the top row of the
diagram is also an equaliser diagram. Now, using that $t_{a
\sqcup_{a_0} \!a}$ is a monomorphism, one can apply Lemma
\ref{2.2.1} to conclude that the square in the diagram is a
pullback showing (c)$\Ra$(b).
\end{proof}

\section{Pairings of functors}

 Generalising the results sketched in the introduction
 we define the notion of pairings of functors on arbitrary categories.

\begin{thm}\label{pair-func}{\bf Pairing of functors.} \em
For any functors $T:\A\to \B$ and $G:\B\to \A$, there is a bijection
(e.g. \cite[2.1]{Par}) between (the class of) natural
transformations between functors
$\A^{op}\times \B\to \Set$,
$$\beta_{a,b}:\A(a, G(b)) \to \A(T(a),b), \quad a\in \A,\, b\in \B,$$
and natural transformations $\sigma: TG\to \id_\B$, with
  $\sigma_a:= \beta_{G(a),a}(\id_{G(a)}): TG(a)\to a$ and
$$\beta_{a,b}:
  a \stackrel{f}\to G(b)\;\longmapsto \; T(a) \stackrel{T(f)}\to TG(b)
    \stackrel{\sigma_b}\to b .$$

We call $(T,G,\sigma)$ a {\em pairing} between the functors $T$ and $G$
and name it a {\em rational pairing} provided the $\beta_{a,b}$ are
monomorphisms for all $a\in \A$ and $b\in \B$.

Clearly, if all the $\beta_{a,b}$ are isomorphisms, then we have
an {\em adjoint pairing}, that is, the functor $G$ is right
adjoint to $T$ and $\sigma$ is just the counit of the adjunction.
Thus rational pairings generalise adjointness.
We mention that, given a pairing $(T,G,\sigma)$,
 Medvedev \cite{Med} calls $T$ a {\em left semiadjoint} to $G$
provided there is a natural transformation $\varphi:\id_A\to GT$ such that
$\sigma T \circ T\varphi = \id_T$. This means that
$\beta_{-,-}$ is a bifunctorial coretraction (dual of \cite[Proposition 1]{Med})
in which case $\sigma$ is a rational pairing.
In case all $\beta_{a,b}$ are epimorphisms, the functor $G$ is said to be
{\em a weak right adjoint} to $T$ in Kainen \cite{Kain}.
 For more about weakened
forms of adjointness we refer to B\"orger and Tholen \cite{BoeTho}.
\end{thm}

Similar to the condition on an adjoint pair of a monad and a comonad
(see \cite{E}) we define

\begin{thm}\label{pair-mon}{\bf Pairing of monads and comonads.} \em
 A \emph{pairing } $\pp=(\bT, {\bG},\sigma)$
between a {\em monad} $\bT=(T, m, e)$ and a
{\em comonad} ${\bG}=(G, \delta, \ve)$ on a
category $\A$ is a pairing $\sigma: TG \to \id$ between the functors $T$ and $G$
inducing - for $a,b\in \A$ - commutativity of the diagrams
\begin{equation}\label{E.5}
\xymatrix{ & \A(a, b) &\\
\A(a, G(b)) \ar[d]_{\A(a,\delta_b)}\ar[rr]^{\beta^\pp_{a,b}}
\ar[ur]^{\A(a,\ve_b)}
&& \A(T(a),b)\ar[d]^{\A(m_a,b)}\ar[lu]_{\A(e_a,\,b)}\\
\A(a, G^2(b))\ar[r]_{\beta^\pp_{a, G(b)}}& \A(T(a), G(b))
\ar[r]_-{\beta^\pp_{T(a), b}}& \A(T^2(a), b),}
\end{equation}
where
\begin{equation}\label{beta}
 \beta^\pp_{a, b}: \A(a, G(b)) \to \A(T(a),b),\quad
 f: a \to G(b)\;\mapsto \;\sigma_b \cdot T(f): T(a) \to b.
\end{equation}
 The pairing $\pp=(\bT,{\bG},\sigma)$ is said to be \emph{rational}
 if $\beta^\pp_{a,b}$ is injective for any $a,b\in \A$.

By the Yoneda Lemma, commutativity of the diagrams in (\ref{E.5})
correspond to commutativity of the diagrams
\begin{equation}\label{E.7}
\xymatrix{ G \ar[r]^-{e G}\ar[dr]_{\ve}& TG\ar[d]^{\sigma}\\ &\id\ ,}
\qquad
\xymatrix{T^2 G \ar[r]^-{T^2\delta}\ar[d]_{mG}& T^2G^2
\ar[r]^-{T\sigma G} & TG \ar[d]^{\sigma}\\
TG \ar[rr]_-{\sigma}&& \id \,  . }
\end{equation}
\end{thm}

\begin{thm}\label{Example 1}{\bf Pairings and morphisms.}
Let $\pp=(\bT, \bG, \sigma)$ be a pairing, $\bT'=(T',m',e')$ any monad,
and $t : \bT' \to \bT$ a monad morphism.
\begin{rlist}
\item The triple $\pp'=(\bT', \bG,
\sigma':=\sigma \cdot tG)$ is also a pairing.
\item If $\pp$ is rational, then $\pp'$ is also rational
provided the natural transformation $t$ is a componentwise
epimorphism.
\end{rlist}
\end{thm}
\begin{proof} (i). The diagram
$$\xymatrix{ G \ar[dr]|{eG}
\ar[d]_{\varepsilon}\ar[r]^-{e'G}& T'G \ar[d]^{tG}\\
\id & TG \ar[l]^{\sigma}}$$
commutes since $t$ is a monad morphism (thus $t \cdot e'=e $) and
because of the commutativity of the triangle in the diagram
(\ref{E.7}.) Thus $\sigma' \cdot  e'G=\sigma \cdot tG \cdot
e'G=\sigma \cdot eG=\varepsilon$.

Consider now the diagram
$$\xymatrix{& T'T'G\ar[r]^-{T'T'\delta}\ar[d]^{ttG}\ar[dl]_{m'G}&
T'T' GG \ar[dr]_{ttGG} \ar[rr]^-{T'tGG}&& T'TGG \ar@{}[rd]|{(4)}
 \ar[dl]^{tTGG}\ar[rr]^-{T'\sigma G}& & T'G \ar[dl]^{tG}\\
T'G \ar[rd]_{tG} \ar@{}[r]|{(1)}& TTG \ar[d]^{mG}
\ar@{}[ru]|{(2)} \ar[rr]_{TT\delta}&& TTGG
\ar@{}[u]|{(3)} \ar[rr]_{T\sigma G}&&TG \ar[dl]^{\sigma}\\
& TG \ar@{}[rru]|{(5)}\ar[rrr]_{\sigma}&&& \id \,,}$$
in which diagram (1) commutes since $t$ is a morphism of monads,
      diagrams (2), (3) and (4) commute by naturality of composition, and
    diagram (5) commutes by commutativity of the rectangle in (\ref{E.7}).
  It then follows that
$$\sigma'\cdot T\sigma' G\cdot T'T'\delta=\sigma \cdot tG \cdot T'\sigma G
\cdot T'tGG \cdot T'T'\delta=\sigma \cdot tG \cdot
m'G=\sigma' \cdot m'G,$$
proving that the triple $\pp'=(\bT', \bG, \sigma'=\sigma \cdot tG)$
is a pairing.
\smallskip

(ii). It is easy to check that the composite
$$\A(a, G(b))\xrightarrow{\beta^{\pp}_{a,b}}\A(T(a),
b)\xrightarrow{\A(t_a,b)}\A(T'(a),b)$$ takes   $f: a \to
G(b)$ to $\sigma'_b \cdot T'(f): T'(a) \to b$ and
thus $-$ since $T(f) \cdot t_a=t_{G(a)}\cdot T'(f)$ by naturality of $t$ $-$ $\beta^{\pp'}_{a,b}=\A(t_a,b)\cdot
\beta^{\pp}_{a,b}$. If $t$ is a componentwise
epimorphism, then $t_a$ is an epimorphism, and then the map
$\A(t_a,b)$ is injective. It follows that
$\beta^{\pp'}_{a,b}$ is also injective provided that
$\beta^{\pp}_{a,b}$ is injective (i.e. the pairing
$\pp=(\bT, \bG, \sigma)$ is rational).
\end{proof}

Dually, one has

\begin{proposition}\label{Example 2}Let $\pp=(\bT, \bG, \sigma)$ be a pairing, $\bG'=(G',\delta',\varepsilon')$ any comonad,
and $t : \bG' \to \bG$ a comonad morphism.
\begin{rlist}
\item The triple $\pp'=(\bT, \bG',
\sigma':=\sigma \cdot Tt)$ is also a pairing.
\item If $\pp$ is rational, then $\pp'$ is also rational
provided the natural transformation $t$ is a componentwise
monomorphism.
\end{rlist}
\end{proposition}

\begin{thm}\label{func-ind}{\bf Functors induced by pairings.}
  Let $\pp=(\bT, {\bG}, \sigma)$ be a pairing on a category $\A$  with
$\beta^\pp_{a,\, b}: \A(a, G(b)) \to \A(T(a),b)$ (see (\ref{beta}).
\begin{zlist}
\item If $(a, \theta_a) \in \A^{\rG}$, then $(a,
\beta^\pp_{a,a}(\theta_a))= (a, \sigma_a \cdot
T(\theta_a)) \in \A_{\rT}$.
\item The assignments
$\begin{array}[t]{rcl}&
(a, \theta_a) &\longmapsto \; (a, \sigma_a \cdot T(\theta_a)), \\
 &f:a\to b &\longmapsto \; f:a\to b,
\end{array}$

yield a conservative functor $\Phi^\pp:\A^{\rG} \to
\A_{\rT}$ inducing a commutative diagram
\begin{equation}\label{d.1.7}
\xymatrix{ \A^{\rG} \ar[rr]^{\Phi^\pp} \ar[dr]_{U^{\rG}} && \A_{\rT} \ar[dl]^{U_{\rT}}\\
&\A\, .& }
\end{equation}
\end{zlist}
\end{thm}

\begin{proof} (1) We have to show that the diagrams
$$
\xymatrix{a\ar[r]^{e_a} \ar@{=}[dr]& T(a)\ar[d]^{\beta^\pp_{a,a}(\theta_a)}\\
& a,} \quad   \xymatrix{ T^2(a) \ar[r]^{m_a}
\ar[d]_{T(\beta^\pp_{a,a}(\theta_a))} &
T(a) \ar[d]^{\beta^\pp_{a,a}(\theta_a)}&\\
T(a) \ar[r]_{\beta^\pp_{a,a}(\theta_a)} & a} $$
are commutative. In the  diagram
$$
\xymatrix{ a \ar[r]^{\theta_a} \ar[d]_{e_a}& G(a) \ar[d]_{e_{G(a)}}
\ar[rd]^{\ve_a}&&\\
T(a) \ar[r]_{T(\theta_a)} & TG(a) \ar[r]_{\sigma_a}&a} $$
the square commutes by naturality of $e: \id \to T$, while the
triangle commutes by (\ref{E.7}). Thus
$\beta^\pp_{a,a}(\theta_a) \cdot e_a= \sigma_a
\cdot T(\theta_a) \cdot e_a=\ve_a \cdot
\theta_a.$ But $\ve_a \cdot \theta_a=\id_a$ since $(a,
\theta_a) \in \A^{\rG}$, implying that
$\beta^\pp_{a,a}(\theta_a) \cdot e_a=\id_a$. This
shows that the left hand diagram commutes.

Since $\beta^\pp_{a,a}(\theta_a)= \sigma_a \cdot
T(\theta_a)$, and $T(\theta_a)\cdot
T(\sigma_a)=T(\sigma_{G(a)}) \cdot T^2G(\theta_a)$ by naturality
of $\sigma$, the right hand diagram can be rewritten as
$$
\xymatrix{T^2(a) \ar[r]^{m_a} \ar[d]_{T^2(\theta_a)}&
T(a)
\ar[r]^{T(\theta_a)} & TG(a) \ar[rd]^{\sigma_a}&&\\
T^2G(a) \ar[r]_{T^2G(\theta_a)} & T^2G^2(a)
\ar[r]_-{T(\sigma_{G(a)})} & TG(a) \ar[r]_-{\sigma_a}& a \, .}$$
It is easy to see that $\sigma_a \cdot T(\theta_a)\cdot m_a
=(\A(m_a, a)\cdot \beta^\pp_{a,a})(\theta_a)$ and
it follows from the commutativity of the bottom diagram in
(\ref{E.5}) that
$$\begin{array}{rl}
\sigma_a \cdot T(\theta_a)\cdot m_a
  &=(\beta^\pp_{T(a), a} \cdot \beta^\pp_{a, G(a)}
  \cdot \A(a,\delta_a))(\theta_a)\\[+1mm]
& =\sigma_a \cdot T(\sigma_{G(a)})\cdot T^2 (\delta_a) \cdot T^2(\theta_a).
\end{array}$$
Recalling that $\delta_a \cdot \theta_a= G(\theta_a) \cdot
\theta_a$ since $(a, \theta_a) \in \A^{\rG}$, we get
$$\sigma_a\cdot T(\theta_a)\cdot m_a=\sigma_a \cdot
T(\sigma_{G(a)})\cdot T^2G(\theta_a)\cdot T^2(\theta_a),$$
proving that the right hand diagram is commutative. Thus
$(a,\beta^\pp_{a,a}(\theta_a)) \in \A_{\rT}.$
\smallskip

(2) By (1), it suffices to show
that if $f: (a, \theta_a) \to (b, \theta_b)$ is a morphism in
$\A^{\rG}$, then $f$ is a morphism in $\A_{\rT}$ from the
$\bT$-module $(a, \sigma_a \cdot T(\theta_a))$ to the
$\bT$-module $(b, \sigma_b \cdot T(\theta_b))$. To say that
$f: a \to b$ is a morphism in $\A_{\rT}$ is to say that the
outer diagram of
$$
\xymatrix{T(a) \ar[r]^{T(\theta_a)\;} \ar[d]_{T(f)} & TG(a)
\ar[r]^{\quad\sigma_a} \ar[d]_{TG(f)}& a \ar[d]^{f}\\
T(b) \ar[r]_{T(\theta_b)\;} & TG(b) \ar[r]_{\quad\sigma_b} & b}$$ is
commutative, which is indeed the case since the left square
commutes because $f$ is a morphism in $\A^{\rG}$, while the right
square commutes by naturality of $\sigma$. Clearly
$\Phi^\pp$ is conservative.
\smallskip

The commutativity of the diagram of functors is obvious.
\end{proof}

\begin{thm}\label{func-prop}{\bf Properties of the functor $\Phi^\pp$.}
Let $\pp=(\bT, {\bG}, \sigma)$ be a pairing on a category $\A$
with induced functor  $\Phi^\pp:\A^{\rG} \to \A_{\rT}$ (see
\ref{func-ind}).
\begin{zlist}
\item The functor $\Phi^\pp$ is comonadic if and only if
      it has a right adjoint.
\item Let $\C$ be a small category such that any functor $\C \to \A$
      has a colimit that is preserved by $T$. Then the functor
      $\Phi^\pp$ preserves the colimit of any functor $\C \to \A^{\rG}$.
\item When $\A$ admits and $T$ preserves all small colimits, the
      category $\A^{\rG}$ has and the functor $\Phi^\pp$
      preserves all small colimits.
\item Let $\A$ be a locally presentable category, suppose that $G$
      preserves filtered colimits and $T$ preserves all small colimits.
      Then $\Phi^\pp$ is comonadic.
\end{zlist}
\end{thm}

\begin{proof}  (1) Since
  $U^G$ creates equalisers of $G$-split pairs,
  $U^T$ creates all equalisers that exist in $\A$, and
  $U_T \Phi^\pp =U_G$,
$\A^G$ admits equalisers of $\Phi^\pp$-split
pairs and $\Phi^\pp$ preserves them. The result now
follows from \ref{func-ind}(2).
\smallskip

(2) Since $T$ preserves the colimit of any functor $\C \to \A$,
the category $\A_{\rT}$ admits and the functor $U_{\rT}$ preserves
the colimit of any functor $\C \to \A_{\rT}$ (see  \cite{Bo}).
Now, if $F: \C \to \A^{\rG}$ is an arbitrary functor, then the
composition $\Phi^\pp \cdot F: \C \to \A_{\rT}$ has a
colimit in $\A_{\rT}$. Since $U_{\rT} \cdot \Phi^\pp
=U^{\rG}$ and $U^{\rG}$ preserves all colimits that exist in
$\A^{\rG}$, it follows that the functor $U_{\rT}$ preserves the
colimit of the composite $\Phi^\pp \cdot F$. Now the
assertion follows from the fact that an arbitrary conservative
functor reflects such colimits as it preserves.
\smallskip

(3) is a corollary of (2).
\smallskip

(4) Note first that Ad\'amek and Rosick\'y proved in
\cite{AR} that the Eilenberg-Moore category with respect to a
filtered-colimit preserving monad on a locally presentable category
is locally presentable.
This proof can be adopted to show that if $G$ preserves
filtered colimits and $\A$ is locally presentable, then
$\A^{\rG}$ is also locally presentable. Therefore
$\A^{\rG}$ is finitely complete, cocomplete, co-wellpowered
and has a small set of generators (see \cite{AR}). Since the
functor $\Phi^\pp$ preserves all small colimits by (3), 
it follows from the (dual of the) Special Adjoint Functor Theorem
(see \cite{M}) that $\Phi^\pp$ admits a right adjoint
functor. Combining this with (1.ii) gives that the functor
$\Phi^\pp$ is comonadic.
\end{proof}

\begin{thm} \label{rat-pair}{\bf Properties of rational pairings.}
Let $\pp=(\bT, {\bG},\sigma)$ be rational pairing on $\A$.
\begin{zlist}
\item The functor $\Phi^\pp: \A^{\rG} \to \A_{\rT}$ is
full and faithful. \item  The functor $G$ preserves monomorphisms.
\item If $\A$ is abelian and $G$ preserves cokernels, then $G$ is
left exact.
\end{zlist}
\end{thm}
\begin{proof}
(1) Obviously, $\Phi^\pp$ is faithful. Let
$(a,\theta_a), (b, \theta_b) \in \A^{\rG}$ and let
$$f:\Phi^\pp(a,\theta_a)=(a, \sigma_a \cdot T(\theta_a)) \to (b, \sigma_b
\cdot T(\theta_b))=\Phi^\pp(b, \theta_b)$$ be a morphism
in $\A_{\rT}$. We have to show that the diagram
\begin{equation}\label{E.11}
\xymatrix{a \ar[r]^f \ar[d]_{\theta_a} & b \ar[d]^{\theta_b}\\
G(a) \ar[r]_{G(f)}& G(b)}
\end{equation} commutes.
Since $f$ is a morphism in $\A_{\rT}$, the diagram
\begin{equation}\label{E.12}
\xymatrix{ T(a) \ar[r]^{T(\theta_a)} \ar[d]_{T(f)}& TG(a)
\ar[r]^{\sigma_a}& a \ar[d]^{f}\\
T(b) \ar[r]_{T(\theta_b)}& TG(b) \ar[r]_{\sigma_b}& b}
\end{equation} commutes and we have
$$\begin{array}{rl}
\beta^\pp_{a,b}(\theta_b \cdot f)& =\;
\sigma_b \cdot T(\theta_b \cdot f)=\sigma_b \cdot T(\theta_b) \cdot T(f) \\
{}_\text{by (\ref{E.12})} & =f \cdot \sigma_a \cdot T(\theta_a) \\
{}_\text{$\sigma$ is a natural}&= \sigma_b \cdot TG(f) \cdot
T(\theta_a)=\sigma_b \cdot T(G(f) \cdot
\theta_a)=\beta^\pp_{a,b}(G(f) \cdot \theta_a).
\end{array}$$
Since $\beta^\pp_{a,b}$ is injective by assumption, it
follows that $G(f) \cdot \theta_a=\theta_b \cdot f$, that is,
(\ref{E.11}) commutes.
\smallskip

(2) Let $f: a \to b$ be a monomorphism in $\A$. Then for
any $x \in \A$, the map
$$\A(T(x), f): \A(T(x), a) \to \A(T(x), b)$$
is injective. Considering the commutative diagram
$$\xymatrix{\A(x, G(a)) \ar[rr]^{\beta^\pp_{x, a}} \ar[d]_{\A(x, G(f))}
&& \A(T(x), a) \ar[d]^{\A(T(x), f)}\\
\A(x, G(b)) \ar[rr]_{\beta^\pp_{x, b}}&& \A(T(x),
b)\,,}$$ one sees that the map $\A(x, G(f))$ is injective for all
$x \in \A$, proving that $G(f)$ is a monomorphism in $\A$.
\smallskip

(3) is a consequence of (2).
\end{proof}

\begin{thm} \label{commorphism} {\bf  Comonad morphisms.} \em
Recall (e.g. \cite{TTT}) that a morphism of comonads $t: {\bG}\to
{\bG}'$ induces a functor
$$t_* : \A^{G} \to \A^{G'},\quad
(a,\theta_a)\mapsto (a, t_a \cdot \theta_a).$$

The passage $t \to t_*$ yields a bijection between comonad
morphisms ${\bG} \to {\bG}'$ and functors $V: \A^{\rG} \to
\A^{\rG'}$ with $U^{G'} V = U^G$.

If $V:\A^{\rG} \to \A^{\rG'}$
is such a functor, then the image of any
cofree ${\bG}$ -comodule $(G(a), \delta_a)$
 under $V$ has the form $(G(a), s_a)$ for some $s_a: G(a)\to G'G(a)$.
Then the collection $\{s_a\,| \, a \in \A\}$ constitute a natural
transformation $s: G \to G'G$ such that $G'\varepsilon \cdot s
:G \to G'$ is a comonad morphism.
\end{thm}

\begin{thm} \label{adj-mon} {\bf Right adjoint for $T$.} \em
Let $\pp=(\bT, {\bG},\sigma)$ be a pairing and suppose
that there exists a comonad $\bT^\di=(T^\di,\delta^\di ,
\varepsilon^\di)$ that is right adjoint to the monad $\bT$
with unit $\overline{\eta}: 1 \to T^\di T$. In this situation,
 the functor $K_{T,T^\di}: \A_{\rT} \to \A^{\rT^\di}$ that takes any
$(a, h_a)\in \A_{\rT}$ to $(a, T^\di(h_a)\cdot
\overline{\eta}_a)$, is an isomorphism of categories and
$U^{\rT^\di}K_{T,T^\di}=U_{\rT}$ (e.g. \cite{MW}). Since
$U_{\rT}\Phi^\pp =U^{\rG}$, one gets the commutative
diagram
$$\xymatrix{ \A^{\rG} \ar[rr]^{\Phi^\pp }
\ar[drr]_{U^{\rG}}&&\A_{\rT} \ar[rr]^{K_{T,T^\di}} &&
\A^{\rT^\di} \ar[lld]^{U^{\rT^\di}}\\
&& \A &.} $$ It follows that there is a morphism of comonads $t:
{\bG} \to \bT^\di$ with $K_{T,T^\di} \Phi^\pp =t_*.$
\end{thm}

\begin{lemma}\label{L.1.14} In the situation given in \ref{adj-mon},
$t$ is the composite $$\xymatrix{G \ar[r]^-{\overline{\eta}G} & T^\di
TG \ar[r]^-{T^\di \sigma}& T^\di\,.}$$
\end{lemma}
\begin{proof}Since
\begin{rlist}
 \item[$\bullet$] for any $(a, \theta_a)\in \A^{\rG}$, $\Phi^\pp (a, \theta_a)=(a,
\sigma_a \cdot T(\theta_a))$, and
 \item[$\bullet$] for any $(a,h_a)\in \A_{\rT}$,
  $K_{T,T^\di}(a,h_a)=(a,T^\di(h_a)\cdot\overline{\eta}_a)$,
\end{rlist}
it follows that for any cofree ${\bG}$-comodule
$(G(a), \delta_a)$,
 $$K_{T,T^\di}\Phi^\pp (G(a),\delta_a )=(G(a), T^\di(\sigma_{G(a)})\cdot T^\di
T(\delta_a)\cdot \overline{\eta}_{G(a)}),$$
thus
$$t_a=T^\di(\ve_a)\cdot T^\di(\sigma_{G(a)})\cdot
T^\di T(\delta_a)\cdot \overline{\eta}_{G(a)}.$$
But since
$$T^\di(\ve_a)\cdot
T^\di(\sigma_{G(a)})=T^\di(\sigma_a)\cdot T^\di
TG(\ve_a)$$ by naturality of $\sigma$ and
$G(\ve_a) \cdot \delta_a=\id_{G(a)} $, one has
$$t_a=T^\di(\sigma_a)\cdot T^\di TG(\ve_a) \cdot
T^\di T(\delta_a)\cdot
\overline{\eta}_{G(a)}=T^\di(\sigma_a) \cdot
\overline{\eta}_{G(a)}.$$
\end{proof}

\begin{thm}\label{right-ad}{\bf Right adjoint functor of $t_*$.} \em
In the setting of \ref{adj-mon},
suppose that the category $\A^G$ admits equalisers. Then it is
well-known (e.g. \cite{D}) that for any comonad
morphism $t: G \to T^\di$, the functor $t_* : \A^G \to \A^{T^\di}$
admits a right adjoint $t^* :\A^{T^\di} \to \A^G$, which can be
calculated as follows. Recall from \ref{mon-comon} that, for any comonad $G$, we denote by $\eta^G$ and $\ve^G$ the unit
and counit of the adjoint pair $(U^G, \phi^G)$.
 Writing $\alpha_t$ for the composite
$$t_* \phi^G \xrightarrow{\eta^{T^\di}t_* \phi^G} \phi^{T^\di}U^{T^\di}t_*\phi^G=\phi^{T^\di}U^G\phi^G
\xrightarrow{\phi^{T^\di} \varepsilon^G}\phi^{T^\di},$$ and
$\beta_t$ for the composite
$$\phi^G U^{T^\di} \xrightarrow{\eta^G \phi^GU^{T^\di}}\phi^G U^G \phi^G U^{T^\di}=\phi^G U^{T^\di}t_*\phi^G U^{T^\di}
\xrightarrow{\phi^GU^{T^\di} \alpha_t
U^{T^\di}}\phi^GU^{T^\di}\phi^{T^\di}U^{T^\di},$$ then $t^*$ is
the equaliser (we assumed that $\A^G$ has equalisers)
$$ \xymatrix{t^* \ar[r]^-{i_t}& \phi^G U^{T^\di}
 \ar@<-0.5ex>[rr]_-{\beta_t}
\ar@<0.5ex>[rr]^-{\phi^GU^{T^\di} \eta^{T^\di}} &&
\phi^GU^{T^\di}\phi^{T^\di}U^{T^\di}.}$$

Note that the counit $\varepsilon_t : t_*t^* \to \id$ of the
adjunction $t_* \dashv t^*$ is the unique natural transformation
that makes the square in the following diagram commute,
\begin{equation}\label{d.1.10}
\xymatrix{t_*t^* \ar@{.}[d]_{\varepsilon_t}\ar[r]^-{t_*i_t}&
t_*\phi^G U^{T^\di}
 \ar[d]_{\alpha_t U^{T^\di}}\ar@<-0.5ex>[rr]_-{t_*\beta_t}
\ar@<0.5ex>[rr]^-{t_*\phi^GU^{T^\di} \eta^{T^\di}} &&
t_*\phi^GU^{T^\di}\phi^{T^\di}U^{T^\di}  \ar[d]^{\alpha_t U^{T^\di}\phi^{T^\di}U^{T^\di}}\\
\id \ar[r]_-{\eta^{T^\di}}& \phi^{T^\di} U^{T^\di}
 \ar@<-0.5ex>[rr]_-{\eta^{T^\di}\phi^{T^\di}U^{T^\di}}
\ar@<0.5ex>[rr]^-{\phi^{T^\di}U^{T^\di} \eta^{T^\di}} &&
\phi^{T^\di}U^{T^\di}\phi^{T^\di}U^{T^\di}.}
\end{equation}

It is not hard to see that for any $a \in \A$, the $a$-component of
the natural transformation $\alpha_t$ is just the morphism $t_a
:G(a) \to T^\di(a)$ seen as a morphism
$$t_* \phi^G(a)=(G(a), t_{G(a)}\cdot \delta_a) \to \phi^{T^\di}(a)=(T^\di(a),
 \delta^\di_a)$$
 in $\A^{T^\di}$. Indeed, since for any $a \in
A$, $(\varepsilon^G)_a=\varepsilon_a$, while for any
$(a,\nu_a) \in \A^{T^\di}$, $(\eta^{T^\di})_{(a,\nu_a)}=\nu_a$,
$(\alpha_t)_a$ is the composite $T^\di(\varepsilon_a)\cdot
t_{G(a)}\cdot \delta_a$. Considering now the diagram
$$\xymatrix{G(a) \ar[r]^-{\delta_a} \ar@{=}[rd] & GG(a)
\ar[r]^-{t_{G(a)}}\ar[d]^{G(\ve_a)}& T^\di G(a)
\ar[d]^{T^\di(\ve_a)}\\
& G(a) \ar[r]_{t_a}& T^\di(a)\,,}
$$ in which the square commutes by naturality of composition,
while the triangle commutes by the definition of a comonad, one sees
that $(\alpha_t)_a=t_a$.
\end{thm}

\begin{theorem}\label{T.1.15} Let $\pp =(\bT, {\bG},\sigma)$
 be a pairing on $\A$. Suppose that $\A^G$ admits equalisers
and that the monad $\bT$ has a right adjoint comonad
$\bT^\di$. Then
\begin{zlist}
\item the functor $\Phi^\pp : \A^{\rG} \to \A_{\rT}$ (and hence
also $t_* : \A^G \to \A^{T^\di}$) is comonadic;
\item  if the pairing $\pp$ is rational, then $\A^{\rG}$ is equivalent to
a reflective subcategory of $\A_{\rT}.$
\end{zlist}
\end{theorem}

\begin{proof} By \ref{func-prop}(1), $\Phi^\pp $ is comonadic if and only if it has a right adjoint.

(1) Since the category $\A^G$ admits equalisers, the functor
$t_* : \A^{\rG} \to \A^{\rT^\di}$ has a right adjoint
$t^* : \A^{T^*} \to \A^{\rG}$ (by \ref{right-ad}). Then evidently
 the functor $t^*K_{T,G}$ is right adjoint to the
functor $\Phi^\pp $. Since $K_{T,T^\di}\Phi^\pp=t_*$,
it is clear that $t_*$ is also comonadic. This completes
the proof of the first part.

(2) If $\pp$ is rational,
$\Phi^\pp $ is full and faithful by \ref{rat-pair}(1), 
and when $\Phi^\pp $ has a right adjoint, the unit of the
adjunction is a componentwise isomorphism (see \cite{M}).
\end{proof}

Given two functors $F, F' : \A \to \B$, we write $\Nat(F, F')$ for the
collection of all natural transformations from $F$ to $F'$.
As a consequence of the Yoneda Lemma recall:

\begin{lemma}\label{L.1.17} Let $\bT=(T,m,e)$ be a monad
on the category $\A$ with right adjoint comonad
$\bT^\di=(T^\di, \delta^\di,\varepsilon^\di)$, unit
$\overline{\eta}: \id \to T^\di T$ and counit $\overline{\ve}:TT^\di\to \id$.
Then for any endofunctor $G : \A \to \A$, there is a bijection
$$\chi: \Nat(TG, \id)\to \Nat(G, T^\di), \quad
 TG \stackrel{\sigma}\to \id \; \longmapsto \;
 G\stackrel{\overline{\eta}G} \to T^\di TG \stackrel{\sigma}\to T^\di,
$$
 with the inverse given by the assignment
$$G\stackrel{s}\to T^\di\;\longmapsto\;
 TG\stackrel{Ts}\to TT^\di \stackrel{\overline{\ve}}\to \id.$$
\end{lemma}

\begin{proposition}\label{P.1.18} Let $\bT$ be a monad on a
category $\A$ with right adjoint comonad $\bT^\di$ (as in
\ref{L.1.17}). Then, for a comonad ${\bG}=(G,\delta,\ve)$
on $\A$ and a natural transformation $\sigma: TG \to \id$,
the following are equivalent:
\begin{blist}
 \item  The triple $\pp =(\bT, {\bG}, \sigma)$ is a pairing.
 \item   $\chi(\sigma) : G \to T^\di$ is a morphism of comonads.
\end{blist}
\end{proposition}
\begin{proof}  $\text{(a)}\Rightarrow \text{(b)}$ follows from Lemma
\ref{L.1.14}, while $\text{(b)}\Rightarrow \text{(a)}$ follows
from the dual of \ref{Example 1}.
\end{proof}

\begin{proposition}\label{P.1.19}
Let $\bT$ be a monad on $\A$ with right adjoint comonad $\bT^\di$ (as in \ref{L.1.17})
and let ${\bG}=(G,\delta, \varepsilon)$ be a  comonad on $\A$.
\begin{zlist}
 \item There exists a bijection between
\begin{rlist}
    \item natural transformations $\sigma: TG \to \id$ for which the
    triple $\pp =(\bT, {\bG}, \sigma)$ is a pairing;

    \item comonad morphisms ${\bG} \to \bT^\di$;

    \item functors $V: \A^{\rG} \to \A_{\rT}$ such that
    $U_{\rT} V=U^{\rG}$.
\end{rlist}
\item A functor $V: \A^{\rG} \to
\A_{\rT}$ with $U_{\rT} V=U^{\rG}$ is an isomorphism
if and only if there exists an isomorphism of comonads
${\bG} \simeq  \bT^\di$.
\end{zlist}
\end{proposition}

\begin{proof}
(1) By Proposition \ref{P.1.18}, it is enough to show
that there is a bijective correspondence between comonad morphisms
${\bG} \to \bT^\di$ and functors $V: \A^{\rG} \to \A_{\rT}$ such
that $U_{\rT} V=U^{\rG}$. But to give a comonad morphism ${\bG} \to
\bT^\di$ is to give a functor $W: \A^{\rG} \to \A^{\rT^\di}$ with
$U^{\rT^\di}W=U^{\rG}$, which is in turn equivalent - since
$K_{T,T^\di}^{-1}$ is an isomorphism of categories with $U_{\rT}
K_{T,T^\di}^{-1}=U^{\rT^\di}$ (see \ref{adj-mon}) - to giving a
functor $V:\A^{\rG} \to \A_{\rT}$ with $U_{\rT} V=U^{\rG}$.
\smallskip

 (2)  According to (1), to give a
functor $V: \A^{\rG} \to \A_{\rT}$ with $U_{\rT} V=U^{\rG}$ is to give
a comonad morphism $t:{\bG} \to \bT^\di$ such that
$t_*=K_{T,T^\di}V$. It follows - since $K_{T,T^\di}$ is an
isomorphism of categories - that $V$ is an isomorphism of
categories if and only if $t_*$ is, or, equivalently, if $t$ is an
isomorphism of comonads.
\end{proof}

\begin{proposition}\label{P.1.21}
Let $\pp =(\bT, {\bG},\sigma)$ be a pairing on a category
$\A$ and $\bT^\di=(T^\di, \delta^\di ,\ve^\di)$
a comonad right adjoint to $\bT$. Consider
the statements:
\begin{rlist}
    \item the pairing $\pp $ is
    rational;
    \item  $\chi(\sigma): G \to T^\di$ is componentwise a monomorphism;
    \item the functor $\Phi^\pp : \A^{\rG} \to \A_{\rT}$ is
    full and faithful.
\end{rlist}
Then one has the implications $\emph{(i)}
\Leftrightarrow \emph{(ii)}\Rightarrow \emph{(iii)}.$
\end{proposition}

\begin{proof}  $\text{(i)} \Rightarrow\text{(iii)}$ is just
  \ref{rat-pair}(1).

 $\text{(i)}\Leftrightarrow \text{(ii)}$  For
any $a \in \A$, consider the natural transformation
$$\xymatrix{\A(a, G(-))\ar[r]^-{\beta^\pp _{a, -}}&
\A(T(a),- ) \ar[r]^-{\alpha_{a, -}} & \A(a,T^\di(-))},$$ where
$\alpha$ denotes the isomorphism of the adjunction $T\dashv
T^\di$. It is easy to see that the induced natural transformation
$G \to T^\di$ is just $\chi(\sigma)$. Since each component of
$\alpha_{a,-}$ is a bijection, it follows that $\chi(\sigma)$ is a
componentwise monomorphism if and only if each component of
$\beta^\pp _{a, -}$ is monomorphism for all $a \in \A$.
\end{proof}

\section{Rational functors}

Let $\A$ be an arbitrary category admitting pullbacks.

\begin{thm}\label{ass-3} \em
Throughout this section we fix a rational pairing $\pp
=(\bT, {\bG}, \sigma)$ on the category $\A$ with a comonad
$\bT^\di=(T^\di, \delta^\di , \varepsilon^\di )$ right
adjoint to $\bT$, unit $\overline{\eta}: 1 \to T^\di T$ and counit
$\overline{\varepsilon}: TT^\di \to \id$.  Let $t=\chi(\sigma)$ (see
Lemma \ref{L.1.17}).
\end{thm}

For any $(a, \vartheta_a) \in \A^{\bT^\di }$, write $\Upsilon(a,
\vartheta_a)$ for a chosen pullback

\begin{equation}\label{D.16}
\xymatrix{ \Upsilon(a, \vartheta_a) \ar[r]^{p_2} \ar[d]_{p_1} & G(a) \ar[d]^{t_a}\\
a \ar[r]_{\vartheta_a}& T^\di(a)\,.}
\end{equation} Since monomorphisms are stable under pullbacks in any category and
since $t_a$ and $\vartheta_a$ are both monomorphisms, it follows
that $p_1$ and $p_2$ are also monomorphisms.

As a right adjoint functor, $T^\di$ preserves all limits existing
in $\A$ and thus the category $\A^{\bT^\di }$ admits those limits
existing in $\A$. Moreover, the forgetful functor $U^{T^\di }:
\A^{T^\di }\to \A$ creates them; hence these limits
(in particular  pullbacks) can be computed in $\A$.

Now,  $\vartheta_a : a \to T^\di(a)$ is the $(a,
\vartheta_a)$-component of the unit $\eta^{T^\di}: \id \to
\phi^{T^\di}U^{T^\di}$ and thus it can be seen as a morphism in
$\A^{\bT^\di }$ from $(a, \vartheta_a)$ to $(T^\di(a),
(\delta^\di )_a)$, while $t_a : G(a) \to T^\di(a)$ is the
$U^{T^\di}(a, \vartheta_a)$-component of the natural
transformation $\alpha_t: t_*\phi^{T^\di}\to \phi^{T^\di}$ (see
\ref{right-ad}) and thus it can be seen as a morphism in
$\A^{\bT^\di }$ from $t_*(G(a), \delta_a)=(G(a),t_{G(a)}\cdot \delta_a) $
to $(T^\di(a), \delta^\di_a)$.
It follows that the diagram (\ref{D.16}) underlies a pullback
in $\A^{\bT^\di }$. In other words, there exists exactly one
$T^\di$-coalgebra structure $\vartheta_{\Upsilon(a, \vartheta_a)}
:\Upsilon(a, \vartheta_a) \to T^\di(\Upsilon(a, \vartheta_a)) $ on
$\Upsilon(a, \vartheta_a)$ making the diagram
\begin{equation}\label{D.17}
\xymatrix{(\Upsilon(a, \vartheta_a),\vartheta_{\Upsilon(a,
\vartheta_a)}) \ar[rr]^-{p_2} \ar[d]_{p_1}&& (G(a),t_{G(a)}\cdot (\delta_G)_a) \ar[d]^{t_a}\\
(a,\vartheta_a) \ar[rr]_{\vartheta_a} &&
(T^\di(a),(\delta^\di )_a) }\end{equation}
 a pullback in
$\A^{\bT^\di }$. Moreover, since in any functor category,
pullbacks are computed componentwise, it follows that the diagram
(\ref{D.17}) is the $(a, \vartheta_a)$-component of a pullback
diagram in $\A^{T^\di}$,
\begin{equation}\label{D.2.3}
\xymatrix{\Upsilon \ar[r]^-{P_2} \ar[d]_{P_1}& t_* \phi^G
U^{T^\di}\ar[d]^{\alpha_t U^{T^\di}}\\
\id \ar[r]_-{\eta^{T^\di}}& \phi^{T^\di}U^{T^\di} .}
\end{equation}
Since the forgetful functor $U_T : \A_T \to \A$ respects
monomorphisms and $U^{T^\di} \, K_{T,T^\di }=U_T$, it
follows that the forgetful functor $U^{T^\di} :\A^{T^\di} \to \A$
also respects monomorphisms. Thus, the natural transformations
$\alpha_t U^{T^\di}$ and $\eta^{T^\di}$ are both componentwise
monomorphisms and hence so too is the natural transformation $P_1:
\Upsilon \to \id.$

Summing up, we have seen that for any $(a, \vartheta_a) \in \A
^{\bT^\di }$, $\Upsilon(a, \vartheta_a)$ is an object of
$\A^{T^\di }$ yielding an endofunctor
$$\Upsilon: \A^{T^\di } \to \A^{T^\di}, \quad
(a, \vartheta_a) \mapsto \Upsilon(a, \vartheta_a).$$

As we shall see later on, the endofunctor $\Upsilon$ is
- under some assumptions -
the functor-part of an idempotent comonad on $\A^{T^\di }$.

\begin{proposition}\label{P.2.4.}
Under the assumptions from \ref{ass-3},
 suppose that $\A$ admits and $G$ preserves equalisers.
Then the category $\A^G$ also admits
equalisers and the functor $t_* : \A^G \to \A^{T^\di}$  preserves them.
\end{proposition}
\begin{proof} Since
the functor $G$ preserves equalisers, the forgetful
functor $U^G : \A^G \to \A$  creates and preserves equalisers.
Thus $\A^G$ admits equalisers.

Since the forgetful functor $U^{T^\di} : \A^{T^\di} \to \A$ also
creates and preserves equalisers, it  follows from the
commutativity of the diagram
$$
\xymatrix{ \A^G \ar[r]^-{t_*} \ar[rd]_{U^G}
&\A^{T^\di} \ar[d]^{U^{T^\di}}\\
& \A}$$ that the functor $t_* : \A^G \to \A^{T^\di}$ preserves
equalisers.
\end{proof}

Note that when $\A$ is an abelian category, and $G$ preserves
cokernels, then $G$ preserves equalisers by Proposition
\ref{p.2.2} and by \ref{rat-pair}(3).

Note also that it follows from the previous proposition that if
$\A$ admits and $G$ preserves equalisers, then the category $\A^G$
admits equalisers. In view of Proposition \ref{P.2.4.}, it then
follows from \ref{right-ad} that the functor $t_*: \A^G \to
\A^{T^\di}$ has a right adjoint $t^* :\A^{ T^\di} \to \A^G$.

\begin{proposition}\label{P.2.5b} Under the conditions of Proposition \ref{P.2.4.},
 the functor $\Upsilon : \A^{T^\di} \to \A^{T^\di}$ is isomorphic to
the functor part $t_*t^*$ of the $\A^{T^\di}$-comonad generated by
the adjunction $t_* \dashv t^* :\A^{T^\di} \to \A^G.$
\end{proposition}

\begin{proof}Since the functor $G$ preserves equalisers,
the functor $t_* : \A^G \to \A^{T^\di}$ also preserves equalisers
by Proposition \ref{P.2.4.}. Then the top row in the diagram
(\ref{d.1.10}) is an equaliser, and since the bottom row is also
an equaliser and the natural transformation $\alpha_t :
t_*\phi^G \to \phi^{T^\di}$ is a componentwise monomorphism (since
$(\alpha)_a=t_a$ for all $a \in \A$, see \ref{right-ad}), it
follows from Lemma \ref{2.2.1} that the square in the diagram
(\ref{d.1.10}) is a pullback. Comparing this pullback with
$(\ref{D.2.3})$, one sees that $\Upsilon$ is isomorphic to
$t_*t^*$ (and $P_1$ to $\varepsilon_t$).
\end{proof}

Note that in the situation of the previous proposition,
$\varepsilon_t: \Upsilon \to \id$ is componentwise a  monomorphism.
\smallskip

For the next results we will assume that the $\A$ has and $G$ preserves
equalisers. This implies in particular that the category $\A^G$
admits equalisers.

\begin{proposition}\label{P.2.5}
With the data given in \ref{ass-3} assume that $\A$ has and $G$ preserves equalisers. Let $(a, \vartheta_a)$ be an any object of $\A^{\bT^\di}$. Then
\begin{rlist}
\item   $\Upsilon (\Upsilon (a, \vartheta_a)) \simeq \Upsilon (a,
\vartheta_a).$
\item  For every regular $\bT^\di $-subcomodule $(a_0,\vartheta_{a_0})$
of $(a, \vartheta_a),$ the following diagram is a pullback,
$$\xymatrix{\Upsilon
(a_0,\vartheta{a_0})\ar[d]_{(\varepsilon_t)_{(a_0,\vartheta{a_0})}}
\ar[r]^-{\Upsilon(i)}
& \Upsilon (a, \vartheta_{a, \vartheta_a})\ar[d]^{(\varepsilon_t)_{(a, \vartheta_a)}}\\
a_0 \ar[r]_-{i}& a.}$$
\end{rlist}
\end{proposition}
\begin{proof}
(i) As we have observed after Proposition \ref{P.2.5b}, the
natural transformation $\varepsilon_t: \Upsilon \to \id$ is
componentwise monomorphism. Thus   $\Upsilon$ is an idempotent
endofunctor, that is, $\Upsilon (\Upsilon (a, \vartheta_a)) \simeq
\Upsilon (a, \vartheta_a).$

(ii) For any regular monomorphism $i: (a_0, \vartheta_{a_0}) \to (a,
\vartheta_a)$ in $\A^{\bT^\di }$, consider the
commutative diagram
\begin{equation}\label{D.18}
\xymatrix{ \Upsilon(a_0, \vartheta_{a_0})
\ar[dd]_{(\varepsilon_t)_{(a_0,\vartheta_{a_0})}}
\ar[rd]^{(i_t)_{(a_0, \vartheta_{\!a_0})}} \ar[rr]^{\Upsilon(i)}
&& \Upsilon(a, \vartheta_a)
\ar@{.>}[dd]^(.75){(\varepsilon_t)_{(a,\vartheta_a)}}
\ar[rd]^{(i_t)_{(a, \vartheta_a)}}\\& G(a_0)
\ar[dd]^(.75){t_{a_0}} \ar[rr]^(0.35){G(i)} && G(a) \ar[dd]^{t_a}\\
a_0 \ar@{.>}[rr]^(.35){i} \ar[rd]_{\vartheta_{a_0}} && a
\ar@{.>}[rd]^{\vartheta_a}\\& T^\di(a_0) \ar[rr]^{T^\di(i)}&&
T^\di(a)\, .}
\end{equation}
We claim that the square $(\varepsilon_t)_{(a, \vartheta_a)} \cdot
\Upsilon(i)=i \cdot (\varepsilon_t)_{(a_0,\vartheta_{a_0})}$ is a
pullback. Indeed, if $f: x \to a_0$ and $g: x \to \Upsilon(a,
\vartheta_a)$ are morphisms such that $i\cdot
f=(\varepsilon_t)_{(a_0,\vartheta_{a_0})} \cdot g$, then we have
$$T^\di(i) \cdot \vartheta_{a_0} \cdot f=\vartheta_a \cdot i \cdot
f=\vartheta_a \cdot (\varepsilon_t)_{(a_0,\vartheta_{a_0})} \cdot
g=t_a \cdot (i_t)_{(a, \vartheta_a)} \cdot g,$$ and since the
square $t_a \cdot G(i)=T^\di (i) \cdot t_{a_0}$ is a pullback,
there exists a unique morphism $k : x \to G(a_0)$ with $t_{a_0}
\cdot k=\vartheta_{a_0} \cdot f$ and $(i_t)_{(a, \vartheta_a)}
\cdot g= G(i) \cdot k.$ Since the functor $T^\di$, as a right
adjoint functor, preserves regular monomorphisms, the forgetful
functor $U^{T^\di}: \A^{T^\di} \to \A$ also preserves regular
monomorphisms. Thus $i: a_0 \to a$ is a regular monomorphism in
$\A$. Then the square $t_{a_0} \cdot (i_t)_{(a_0,
\vartheta_{\!a_0})} =\vartheta_{a_0}\cdot
(\varepsilon_t)_{(a_0,\vartheta_{a_0})}$ is a pullback by
Proposition \ref{p.2.2}. Therefore, there exists a unique morphism
$k': x \to \Upsilon (a_0, \vartheta_{a_0})$ with $k= (i_t)_{(a_0,
\vartheta_{\!a_0})}\cdot k'$ and $
(\varepsilon_t)_{(a_0,\vartheta_{a_0})}\cdot k'= f.$ To show that
$\Upsilon(i) \cdot k'=g,$ consider the composite
$$ (\varepsilon_t)_{(a,\vartheta_a)} \cdot \Upsilon(i) \cdot k'=i \cdot
(\varepsilon_t)_{(a_0,\vartheta_{a_0})} \cdot k'=i \cdot f=
(\varepsilon_t)_{(a,\vartheta_a)} \cdot g.$$ Since $
(\varepsilon_t)_{(a,\vartheta_a)}$ is a monomorphism, we get $k
\cdot k'=g$. This completes the proof of the fact that the square
$(\varepsilon_t)_{(a, \vartheta_a)} \cdot \Upsilon(i)=i \cdot
(\varepsilon_t)_{(a_0,\vartheta_{a_0})}$ is a pullback.
\end{proof}

\begin{proposition} \label{P.2.9} 
Assume the same conditions as in Proposition \ref{P.2.5}.
The functor $t_* : \A^G \to \A^{T^\di}$ corestricts to an
equivalence between $\A^G$ and the full subcategory of
$\A^{T^\di}$ generated by those ${T^\di}$-coalgebras
$(a,\vartheta_a)$ for which there exists an isomorphism
$\Upsilon(a,\vartheta_a) \simeq (a,\vartheta_a)$. This holds if
and only if there is a (necessarily unique) morphism $x:a \to
G(a)$ with $t_a \cdot x=\vartheta_a$.
\end{proposition}
\begin{proof}Since the category $\A^G$ admits equalisers, it follows from Theorem
\ref{T.1.15} that the functor $t_*$ is comonadic. Thus $\A^G$ is
equivalent to $(\A^{T^\di})^\Upsilon$. But since $\varepsilon_t:
\Upsilon \to \id$ is a componentwise monomorphism, the result
follows from \ref{idempotent}.

Next, $\Upsilon(a,\vartheta_a) \simeq (a,\vartheta_a)$ if and only
if the morphism $p_1$ in the pullback diagram (\ref{D.17}) is an
isomorphism.
In this case the composite $x= p_2 \cdot
(p_1)^{-1}: a \to G(a)$ satisfies the condition of the
proposition. Since $t_a$ is a monomorphism, it is
clear that such an $x$ is unique.

Conversely, suppose that there is a morphism $x:a \to G(a)$ with
$t_a \cdot x=\vartheta_a$. Then it is easy to see -- using that
$t_a$ is a monomorphism-- that the square
$$
\xymatrix{ a \ar[d]_\id \ar[r]^{x} &G(a) \ar[d]^{t_a}\\
a \ar[r]_{\vartheta_a} &  T^\di(a)}$$ is a pullback. It follows
that $\Upsilon(a,\vartheta_a) \simeq (a,\vartheta_a)$.
\end{proof}
\begin{thm}\label{P.2.8}{\bf Rational functor.}
Assume the data from \ref{ass-3} to be given and that $\A$ has and $G$
preserves equalisers. Define the functor
$$\R^\pp: \xymatrix{\A_{\rT}
\ar[r]^{K_{T,T^\di}}&\A^{\bT^\di } \ar[r]^{\Upsilon}& \A^{\bT^\di
} \ar[r]^{K_{T,T^\di}^{-1}}& \A_{\rT}\, .}$$
Then the triple
$(\R^\pp , \varepsilon^\pp , \delta^\pp)$, where $\varepsilon^\pp=K_{T,T^\di}^{-1}\cdot\varepsilon_t \cdot K_{T,T^\di}$
and
$\delta^\pp=K_{T,T^\di}^{-1}\cdot \delta_t \cdot K_{T,T^\di}$,
is a comonad on $\A_T$. Moreover, for any object
$(a, h_a)$ of $\A_{\rT}$,
\begin{rlist}
\item   $\R ^\pp (\R^\pp  (a, h_a))\simeq\R^\pp  (a, h_a);$
\item for any regular $T$-submodule
$(a_0, h_{a_0})$ of $(a, h_a)$, the following diagram is a
pullback,
 $$\xymatrix{\R^\pp(a_0,\vartheta{a_0})
 \ar[d]_{(\varepsilon\!^\pp)_{(a_0,\vartheta{a_0})}}\ar[r]^-{\R\!^\pp(i)}
& \R^\pp (a,\vartheta_{a,\vartheta_a})\ar[d]^{(\ve\!^\pp)_{(a,\vartheta_a)}}\\
a_0 \ar[r]_-{i}& a.}$$
\end{rlist}
\end{thm}
\begin{proof} Observing that  $(\R^\pp ,\varepsilon^\pp, \delta^\pp)$
is the comonad obtained from the comonad $(\Upsilon,
\varepsilon_t, \delta_t)$ along the isomorphism $K_{T,T^\di}^{-1}
: \A^{T^\di} \to \A_T$ (see \cite{S}), the results follow from
Proposition \ref{P.2.5}.
\end{proof}

We call a $T$-module $(a, h_a)$ \emph{rational} if $\R^\pp(a,h_a)\simeq (a,h_a)$
(which is the case if and only if
$(\varepsilon^\pp)_{(a,h_a)}$ is an isomorphism, see
\ref{idempotent}), and write $\R^\pp(\bT)$ for the
corresponding full subcategory of $\A_{\rT}$. Applying Proposition
\ref{P.2.9} gives:

\begin{proposition}\label{P.2.11}
Under the assumptions of Proposition \ref{P.2.8},
let $(a, h_a) \in \A_T$. Then $(a, h_a) \in \R^\pp (\bT)$
if and only if there exists a (necessarily unique) morphism $x : a
\to G(a)$ inducing commutativity of the diagram
$$
\xymatrix{& G(a) \ar[d]^{t_a}\\
a \ar[ru]^{x} \ar[r]_-{\vartheta_a}&  T^\di(a). }$$
\end{proposition}

Putting together the information obtained so far, we obtain as
main result of this section:

\begin{theorem}\label{T.2.15}
Let $\pp=(\bT,{\bG}, \sigma)$ be a rational pairing on a
category $\A$ with a comonad $\bT^\di=(T^\di, \delta^\di ,
\ve^\di )$ right adjoint to $\bT$. Suppose that $\A$
admits and $G$ preserves equalisers.
\begin{zlist}
 \item $\R^\pp(\bT)$ is a coreflective subcategory of $\A_T$,
       i.e. the inclusion $i_\pp:\R^\pp (\bT)
\to \A_{\rT}$ has a right adjoint ${\rm{rat}}^\pp: \A_{\rT}\to \R^\pp (\bT)$.
\item The idempotent comonad on $\A_T$
generated by the adjunction $i_\pp \dashv\rm{rat}^\pp$
is just the idempotent monad
$(\R^\pp, \varepsilon^\pp)$.
\item The functor
$$\Phi^\pp: \A^{\rG} \to \A_T,\quad (a,\vartheta_a) \mapsto
  (a, \sigma_a \cdot T(\vartheta_a)),$$
corestricts to an equivalence of categories
 $R^\pp:\A^{\rG} \to\R^\pp(\bT)$.
\end{zlist}
\end{theorem}

As a special case, consider for $\A$ the category ${_A\M}$ of left
$A$-modules, $A$ any ring. For an $A$-coring $C$ there is a pairing
$(C^*,C,\ev)$. If this is rational it follows by Theorem \ref{T.2.15}
that the $C$-comodules form a coreflective subcategory of ${_{C^*}\M}$
 (see \ref{pair-cor}).

\section{Pairings in monoidal categories}

We begin by reviewing some standard definitions associated with
monoidal categories.

Let $\cV=(\V, \ot, \II)$ be a monoidal category with tensor
product $\ot$ and unit object $\II$. We will freely appeal to
MacLane's coherence theorem (see \cite{M}); in particular, we
write as if the associativity and unitality isomorphisms were
identities. Thus $X \ot( Y\ot Z)=(X \ot Y)\ot Z$
and $\II \ot X=X \ot \II$ for all $X,Y,X \in \V$. We
sometimes collapse $\ot$ to concatenation, to save space.

Recall that an algebra $\uA $ in $\cV$ (or $\cV$-algebra) is an
object $A$ of $\V$ equipped with a multiplication $m_{A}: A
\ot A \to A$ and a unit $e_{A} : \II \to A$  inducing
commutativity of the diagrams
$$
\xymatrix{ A\ot A \ot A\ar[d]_{A \ot m_{A}}
 \ar[rr]^{m_{A}
\ot A}&& A \ot A
\ar[d]^{m_{A}} \\
A \ot A \ar[rr]_{m_{A}} && A,} \quad \xymatrix{   \II \ot
A \ar@{=}[rd]\ar[r]^{e_{A} \ot A} & A \ot A
\ar[d]_{m_{A}}& A \ot \II \ar[l]_{A \ot e_{A}}
\ar@{=}[ld] \\
 &A . & }
$$

Dually, a coalgebra $\oC$ in $\cV$ (or $\cV$-coalgebra) is an
object $C$ of $\V$ equipped with a comultiplication $\delta_{C}:
 C\to C \ot C$, a counit
$\varepsilon_{C} : C \to \II$ and commutative diagrams
$$
\xymatrix{ C \ot C \ot C && C \ot C
\ar[ll]_{\delta_{C} \ot C}
 \\
C \ot C \ar[u]^{C \ot \delta_{C}} && C
\ar[ll]^{\delta_{C}} \ar[u]_{\delta_{C}},}\quad \xymatrix{  \II
\ot C \ar@{=}[rd] & C \ot C \ar[l]_{\varepsilon_{C}
\ot C}\ar[r]^{C \ot
\varepsilon_{C}} & C \ot \II \ar@{=}[ld] \\
 &C . \ar[u]_{\delta_{C}} & }
$$

\begin{definition}\label{def.pair} \em
A triple $\pp=(\uA , \oC, t)$ consisting of a $\cV$-algebra $\uA
$, a $\cV$-coalgebra $\oC$ and a morphism $t: A\ot C \to \II$,
for which the diagrams
\begin{equation}\label{D.25}
\xymatrix{  A  \ot A \ot C \ar[rr]^-{ A \ot A \ot
\delta_{C}} \ar[d]_{ m_{A}\ot C}&&A \ot A \ot C
\ot C \ar[rr]^-{A  \ot t  \ot C}&& A \ot C
\ar[d]_{t} & \ar[l]_-{e_{A}\ot   C} C
\ar[ld]^{\varepsilon_{C}} \\
A  \ot C \ar[rrrr]_-{t}&&&&\II}
\end{equation} commute, is called a \emph{left pairing}.
\end{definition}

A \emph{left action} of a monoidal category $\cV=(\V, \ot,
\II)$ on a category $\X$ is a functor $$-\diamondsuit-: \V \times
\X \to \X, $$ called the \emph{action of} $\V$ \emph{on} $\X$,
along with invertible natural transformations

$$\alpha_{A,B,X}:(A \ot B)\diamondsuit X \to A \diamondsuit (B \diamondsuit
X)\,\,\text{and}\,\, \lambda_X : \II \diamondsuit X \to  X,$$ called
the \emph{associativity} and \emph{unit} isomorphisms,
respectively, satisfying two coherence axioms (see B\'{e}nabou
\cite{Bn}). Again we write as if $\alpha$ and $\lambda$ were
identities.


\begin{thm}\label{Ex.1}{\bf Example.} \em Recall that if $\mathcal{B}$
is a bicategory (in the sense of B\'{e}nabou \cite{Bn}), then for
any $A \! \in \!\text{Ob}(\mathcal{B})$, the triple
$(\mathcal{B}(A, A), \circ, \id_A)$, where $\circ$ denotes the
horizontal composition operation, is a monoidal category, and
that, for an arbitrary $B \in \mathcal{B}$, there is a canonical
left action of $\mathcal{B}(A,A)$ on $\mathcal{B}(B,A)$, given by
$f\diamondsuit g=f \circ g$ for all $f \in \mathcal{B}(A,A)$ and
all $g \in \mathcal{B}(B,A)$. In particular, since monoidal
categories are nothing but bicategories with exactly one object,
any monoidal category $\cV=(\V, \ot, \II)$ has a canonical
(left) action on the category $\V$, given by $A \diamondsuit B=A
\ot B. $
\end{thm}

\begin{thm}\label{act-pair}{\bf Actions and pairings.} \em
Given a left action $-\diamondsuit-: \V \times \X \to \X $  of a
monoidal category $\cV$ on a category $\X$ and an algebra $\uA
=(A, e_{A}, m_{A})$ in $\cV$, one has a monad  $\bT^\X_{\!\uA }$
on $\V$ defined on any $X \in\cV$ by
\begin{rlist}
\item[$\bullet$] $\bT^\X_{\!\uA }(X)=A \diamondsuit X$,
\item[$\bullet$] $(e_{\bT^\X_{\uA }})_X=e_{A} \diamondsuit X : X=\II
\diamondsuit X \to A \diamondsuit X =\bT^\X_{\!\uA }(X)$,
\item[$\bullet$] $(m_{T^\X_{\uA }})_X =m_{A} \diamondsuit X
:\bT^\X_{\!\uA }(\bT^\X_{\!\uA }(X))=A\diamondsuit (A \diamondsuit
X)=(A \ot A)\diamondsuit X \to A \diamondsuit X=\bT^\X_{\!\uA
}(X),$
\end{rlist}
and we write $_\uA \! \X$ for the Eilenberg-Moore category of
$\bT_{\!\uA }$-algebras. Note that in the case of the canonical
left action of $\cV$ on itself, $_\uA \! \cV$ is just  the
category of (left) $\uA $-modules.

Dually, for a $\cV$-coalgebra $\oC=(C,\varepsilon_C,\delta_C)$, a
comonad ${\bG}_{\X}^{\oC}=(C \diamondsuit -,\, \varepsilon_{C}
\diamondsuit -,\, \delta_{C} \diamondsuit -)$ is defined on $\X$
and one has the corresponding Eilenberg-Moore category $^{\oC}\X$;
 for $\X=\cV$ this is just the category of (left)
$\oC$-comodules.

It is easy to see that if $\pp=(\uA,\oC,t)$ is a left
pairing in $\cV$, then the triple $\pp_\X=(\bT^\X_{\!\uA
}, {\bG}_{\X}^{\oC}, \sigma_t)$, where $\sigma_t$ is the natural
transformation
$$ t \diamondsuit -:\bT_{\uA } \cdot {\bG}^{\oC}=
 A\diamondsuit(C \diamondsuit -)=(A \ot C)\diamondsuit- \to \II\ot - =\id,$$
  is a pairing between the monad $\bT^\X_{\!\uA }$ and comonad
${\bG}_{\X}^{\oC}$.

We say that a left pairing $(\mathbf{A},\textbf{B},\sigma)$ is
$\X$-\emph{rational}, if the corresponding pairing
$\pp_\X$  is rational, i.e. if the map
$$\beta^{\pp_\X}_{X,Y}: \X(X, C \diamondsuit Y) \to \X(A \diamondsuit X,Y)$$
taking $f: X \to C \diamondsuit Y$ to the composite $A
\diamondsuit X \xrightarrow{A \diamondsuit f}A \diamondsuit (C
\diamondsuit Y)={A \ot }\diamondsuit Y\xrightarrow{\sigma
\diamondsuit Y}Y$, is injective.

 We will generally drop the $\X$ from the notations
$\bT^\X_{\!\uA}$, ${\bG}_{\X}^{\oC}$ and $\pp_\X$ when
there is no danger of confusion.
\end{thm}

\begin{thm}\label{closed-cat}{\bf Closed categories.} \em
A monoidal category $\cV=(\V, \ot, \tau, \II)$ is said to be
\emph{right closed} if each functor $- \ot X : \V \to \V$ has
a right adjoint $[X, -]: \V \to \V$. So there is a bijection
\begin{equation}\label{E.24}
\pi_{Y, X, Z}: \V (Y \ot X, Z)\simeq
\V(Y,[X,Z]),\end{equation}
with unit $\eta_X^Y : X \to [Y, X \ot Y]$ and counit $e_Z^Y :
[Y, Z]\ot Y \to Z.$

We write $(-)^*$ for the functor $[-, \II]: \V^{op} \to \V$ that
takes $X \in \V$ to $[X, \II]$ and $f:Y \to X$ to the
morphism $[f, \II]: [X, \II] \to [Y, \II]$ that corresponds under
the bijection (\ref{E.24}) to the composite $$\xymatrix{ [X,
\II]\ot Y \ar[rr]^-{[X, \II]\ot f}&& [X, \II] \ot X
\ar[r]^-{e_\II^X}&\II.}$$

Symmetrically, a monoidal category $\cV=(\V, \ot, \tau, \II)$
is said to be \emph{left closed} if each functor  $X \ot - :
\V \to \V$ has a right adjoint $\{X, -\}: \V \to \V$. We write
$\overline{\eta}^Y_X: X \to [Y, Y \ot X]$ and
$\overline{e}^Y_Z : Y \ot \{Y,Z\} \to Z$ for the unit and
counit of the adjunction $X \ot - \dashv \{X, -\}$. One calls
a monoidal category \emph{closed} when it is both left and right
closed. A typical example is the category of bimodules over a
non-commutative ring R, with $\ot_R$ as $\ot$.
\end{thm}

\begin{thm}\label{act-adj}{\bf Actions with right adjoints.} \em
Suppose now that $-\diamondsuit-: \V \times \X \to \X $ is a left
action of a monoidal category $\cV$ on a category $\X$ and that
$\uA =(A, e_{A}, m_{A})$ is an algebra in $\cV$ such that the
functor $A \diamondsuit-:  \X \to \X $ has a right adjoint $\{A, -
\}_\X:  \X \to \X$ (as it surely is when $\cV=\mathcal{B}(A,A)$
and $\X=\mathcal{B}(B,A)$ for some objects $A, B$ of a
\emph{closed} bicategory $\mathcal{B}$, see Example \ref{Ex.1}.)

Recall (e.g. \cite{E}) that given a monad $\bT=(T, m_{T}, e_{T})$
on a category $\X$ and an endofunctor $G: \X \to \X$ right adjoint
to $T$, there is a unique way to make $G$ into a comonad
${\bG}=(G, \delta_{{G}}, \varepsilon_{{G}})$ such that ${\bG}$ is
right adjoint to the monad $\bT$. Since the functor $\{A, -\}_\X:
\X \to \X$ is right adjoint to the functor $A \diamondsuit -: \X
\to \X$ and since $A \diamondsuit -$ is the functor-part of the
monad $\bT\!_{\uA }$, there is a unique way to make $\{A, -\}_\X$
into a comonad ${\bG}(\uA )=(\{A, -\}_\X, \, \delta_{{\bG}(\uA )},
\, \varepsilon_{{\bG}(\uA )})$ such that the comonad ${\bG}(\uA )$
is right adjoint to the monad $\bT\!_{\uA }$.

Dually, for any coalgebra $\oC=(C, \delta_{C}, \varepsilon_{C})$
in $\cV$ such that the functor  $C \diamondsuit-:  \X \to \X $ has
a right adjoint $\{C, - \}_\X:  \X \to \X$ , there exists a monad
$\bT(\oC)$ whose functor-part is $\{C, -\}_\X$ and which is right
adjoint to the comonad ${\bG^\oC}$.
\end{thm}

\begin{thm}\label{P.3.3}{\bf Pairings with right adjoints.} \em
Consider a left action of a monoidal category $\cV$ on a
category $\X$ and a left pairing $\pp=(\mathbf{A},
\mathbf{C}, \sigma)$ in $\cV$, such that there exist adjunctions
$$A \diamondsuit - \dashv \{A,- \}_\X, \quad - \diamondsuit A \dashv
[A, -]_\X, \quad C \diamondsuit - \dashv \{C,- \}_\X,
 \quad \text{and}\,\,- \diamondsuit C \dashv [C, -]_\X.$$

Since the comonad ${\bG}(\uA )$ is right adjoint to the
monad $\bT_{\uA }$, the following are equivalent
 by Proposition \ref{P.1.21}:
\begin{blist}
  \item the pairing $\pp$ is $\X$-rational;
  \item for every $X \in \X$, the composite
$$ \xymatrix{\alpha^\pp_X :C \diamondsuit X
  \ar[r]^-{\overline{\eta}^A_{C \diamondsuit X}} & \{A, \,
A \diamondsuit (C \diamondsuit X)\}_\X= \{A, \,( A \ot C)
\diamondsuit X\}_\X\ar[rr]^-{\{A, \, t \diamondsuit X\}_\X} &&
\{A, X\}_\X,}$$ where $\overline{\eta}^A_{C \diamondsuit X}$ is
the $C \diamondsuit X$-component of the unit of the adjunction $A
\diamondsuit - \dashv \{A,-\}_\X$, is a monomorphism.
\end{blist}
If (one of) these conditions is satisfied, then the functor
$$\Phi^\pp: {^{\oC}\X} \to {_{\uA }\X}, \;
(X,X \xrightarrow{\theta_X} C \diamondsuit X )\longmapsto(X, A\diamondsuit X \xrightarrow{A \diamondsuit
\theta_X}A \diamondsuit (C \diamondsuit X)=(A \ot C)
\diamondsuit X \xrightarrow{t \diamondsuit X}X)$$
is full and faithful.
\end{thm}

Writing $\R^\pp$ (resp. $\R^\pp(\uA )$) for the functor
$\R^{\pp_\X}$ (resp. the category $\R^{\pp_\X}(\bT_{\uA }))$, one
gets:

\begin{proposition}\label{P.3.4} Under the conditions given in \ref{P.3.3},
assume the category $\X$ to admit equalisers.
If $\pp$ is an $\X$-rational pairing in $\cV$ such that either
\begin{rlist}
 \item the functor $C \diamondsuit -: \X \to \X$
       preserves equalisers, or
 \item $\X$ admits pushouts and the functor $C\diamondsuit -: \X \to \X$
       preserves regular monomorphisms, or
 \item $\X$ admits pushouts and every monomorphism in $\X$ is regular,
\end{rlist}
then
\begin{zlist}
\item the inclusion $i_{\pp}:\R^{\pp}(\uA ) \to {_{\uA }\X}$
  has a right adjoint $\rm{rat}^{\pp}:{_{\uA}\X}\to\R^{\pp}(\uA )$;
\item  the idempotent comonad on ${_{\uA }\X}$ generated by the
adjunction is just $\R^{\pp} : {_{\uA }\X} \to {_{\uA }\X}$; \item
the functor $\Phi^{\pp}: {^{\oC}\X} \to {_{\uA }\X} $ corestricts
to an equivalence $R^{\pp}: \R^{\pp}(\uA ) \to {^{\oC}\X} $.
\end{zlist}
\end{proposition}
\begin{proof} For condition (i), the assertion follows from Theorem \ref{T.2.15}.
We show that (iii) is a particular case of (ii), while (ii) is itself a
particular case of (i). Indeed, since $\pp$ is an
$\X$-rational pairing in $\cV$, it follows from \ref{rat-pair}(2)
that the functor $C \diamondsuit -: \X \to \X$ preserves
monomorphisms, and if every monomorphism in $\X$ is regular, then
$C \diamondsuit -: \X \to \X$ clearly preserves regular
monomorphisms. Next, since the functor $C \diamondsuit -: \X \to
\X$ admits a right adjoint, it preserves pushouts, and then it
follows from Proposition \ref{p.2.2} that $C \diamondsuit -: \X
\to \X$ preserves equalisers. This completes the proof.
\end{proof}

\begin{thm} {\bf Nuclear objects.} \em
We call an object $V \in \V$ is (left) $\X$-\emph{prenuclear}
(resp. $\X$-\emph{nuclear}) if

\begin{zlist}
 \item[$\bullet$] the functor $- \ot V : \V \to \V $ has a right
    adjoint $[ V,-] : \V \to \V, $
 \item[$\bullet$] the functor $V^*\diamondsuit -: \X \to \X$, with $V^*=[V,\II]$,
    has a right adjoint $\{V^*,-\}_\X: \X \to \X$, and
 \item[$\bullet$] the composite
$$ \xymatrix{\alpha_X :V \diamondsuit X \ar[r]^-{(\eta_\X)_{V \diamondsuit X}}
 & \{V^*, \, V^* \diamondsuit (V \diamondsuit X)\}= \{V^*, \,( V^* \ot V)
\diamondsuit X\}\ar[rr]^-{\{V^*, \, e^V_{\II} \diamondsuit X\}} &&
\{V^*, X\},}$$
is a monomorphism (resp. an isomorphism), where
$((\eta)_\X)_{V \diamondsuit X} $ is the $V \diamondsuit
X$-component of the unit $\eta_\X : - \to \{V^*, V^* \diamondsuit
-\}$ of the adjunction $V^* \diamondsuit - \dashv \{V^*,-\}$.
\end{zlist}

Note that the morphism $\alpha_X :V \diamondsuit X \to \{V^*, X\}$
is the transpose of the morphism $e^V_{\II} \diamondsuit X : ( V^*
\ot V) \diamondsuit X \to X$ under the adjunction $V
\diamondsuit- \dashv \{V,-\}_\X$.

\end{thm}

Applying Proposition \ref{P.3.4} and \ref{bicl-coalg}, we get:

\begin{proposition}\label{P.3.7}  Let $\cV$ be a monoidal closed
category and $\oC=(C,\varepsilon_C,\delta_C)$ a $\cV$-coalgebra
with $C$ $\X$-prenuclear, and assume $\X$ to admit equalisers. If
either
\begin{rlist}
    \item [\emph{(i)}] the functor $C \diamondsuit -: \X \to \X$
    preserves equalisers, or
    \item [\emph{(ii)}]$\X$ admits pushouts and the functor $C \diamondsuit -: \X \to \X$
    preserves regular monomorphisms, or
    \item [\emph{(iii)}]$\X$ admits pushouts and every
    monomorphism in $\X$ is regular,
\end{rlist}
then $\R^{\pp(\oC)}(\oC^* ) $ is a full coreflective subcategory
of $_{\oC^*}\X$ and the functor $\Phi^{\pp(\oC)}: {^{\oC}\X} \to
{_{\oC^*}\X} $ corestricts to an equivalence $R^{\pp(\oC)}:
{^{\oC}\X}
 \to \R^{\pp(\oC)}(\oC^* ) $.
\end{proposition}

Specialising the previous result to the case of the left action of
the monoidal category $_{C^*}\!\M_{C^*}$ of $C^*$-bimodules on the
category of left $C^*$-modules, one sees that the equivalence of
the category of comodules $^C\M$ and a full subcategory of
$_{C^*}\M$ for $_AC$ locally projective addressed in
\ref{pair-cor} is a special case of the preceding theorem.
\smallskip

Recall (e.g. from \cite{M}) that a monoidal category
$\cV=(\V, \ot, \II)$ is said to be {\em symmetric} if for all $X, Y
\in \V$, there exists functorial isomorphisms $\tau_{X,Y} : X
\ot Y \to Y \ot X$ obeying certain identities. It is clear
that if $\V$ is closed, then $\{X, -\}\simeq [X,-]$ for all $X \in
\V$.

\begin{thm}\label{symm-pair}{\bf Pairings in symmetric monoidal closed categories.} \em
Let $\cV=(\V, \ot, \II, \tau, [-, -])$ be a symmetric monoidal
closed category and $\pp=(\mathbf{A},\mathbf{C},t)$ a left pairing
in $\cV$. For any $X \in \cV$, we write
$$\gamma_X : A^* \ot X \to [A, X]$$
for the morphism that corresponds under $\pi$ (see (\ref{E.24}))
to the composition
$$\xymatrix{A^* \ot X \ot A \ar[rr]^-{\tau_{A^*,
\, X} \ot A}&& X \ot A^* \ot A \ar[r]^-{e^A_\II}&
X\,.}$$
Since the functor $[A,-]$, as a left adjoint, preserves
colimits, it follows from \cite[Theorem 2.3]{F} that there exists
a unique morphism $\gamma(t): C \to A^*$ such that the diagram
\begin{equation}\label{D.27}
\xymatrix{C \ot X \ar[rr]^-{\gamma(t)\ot X}
\ar[rrd]_-{\alpha_X}&& A^* \ot X \ar[d]^{\gamma_X}\\&&
[A,X]}\end{equation} commutes. Hereby the morphism $\gamma(t)$
corresponds under $\pi$ to the composite $\sigma \cdot
\tau_{C,\,A}$.
\end{thm}

\begin{proposition}\label{P.3.8} Consider the situation given in \ref{symm-pair}.
\begin{zlist}
\item If $\pp$ is a rational pairing, then the morphism $\gamma(t): C \to A^*$ is
pure, that is, for any $X \in \cV$, the morphism $\gamma(t)\ot
X $ is a monomorphism.
\item If $A$ is $\cV$-nuclear, then $\pp$ is rational
      if and only if $\gamma(t): C \to A^*$ is pure. 
\item
$\Phi^\pp: {^{\oC}\cV} \to {_{\uA}\cV}$ is an
isomorphism if and only if $A$ is $\cV$-nuclear and $\gamma(t): C
\to A^*$ is an isomorphism. In this case $C$ is also $\cV$-nuclear.
\end{zlist}
\end{proposition}

\begin{proof} (1) To say that $\pp$ is
a rational pairing is to say that $\alpha_X$ is a monomorphism for
all $X \in \cV$ (see Proposition \ref{P.3.3}). Then it follows
from the commutativity of the diagram (\ref{D.27}) that
$\gamma(t)\ot X $ is also a monomorphism, thus $\gamma(t): C
\to A^*$ is pure.
\smallskip

(2)  follows from the commutativity of diagram (\ref{D.27}).
\smallskip

(3) One direction is clear, so suppose that the functor $\Phi$ is
an isomorphism of categories. Then it follows from Proposition
\ref{P.1.19}(2) that $\alpha_X$ is an isomorphism for all $X \in
\cV$. It implies that the functor $[A, -]$ preserves colimits and
thus the morphism $\gamma_X$ is an isomorphism (see \cite[Theorem
2.3]{F}). Therefore $A$ is $\cV$-nuclear. Since $\alpha_X$ and
$\gamma_X$ are both isomorphisms, it follows from the
commutativity of diagram (\ref{D.27}) that the morphism
$\gamma(t)\ot X$ is also an isomorphism. This clearly implies
that $\gamma(t): C \to A^*$ is an isomorphism.

It is proved in \cite[Corollary 2.2]{R} that if $A$ is
$\cV$-nuclear, then so is $A^*$. 
\end{proof}

\section{Entwinings in monoidal categories}

Recall (for example, from \cite{Me}) that an entwining in a
monoidal category $\cV=(\V, \ot, \II)$ is a triple
$(\uA, \uC,\lambda)$, where $\uA=(A, e_A,
m_A)$ is a $\cV$-algebra, $\uC=(C,\varepsilon_C, \delta_C)$
is a $\cV$-coalgebra, and $\lambda: A \ot C  \to  C \ot A$
is a morphism inducing commutativity of the diagrams
$$
\xymatrix{  C  \ar[d]_{e_A \ot C} \ar[dr]^{C \ot e_A} && A \ot C
\ar[d]_{\lambda} \ar[rd]^-{A \ot \varepsilon_{\uC} }&\\
A \ot C \ar[r]_{\lambda}& C \ot A \,,  &  C \ot A
\ar[r]_{\varepsilon_{\uC} \ot A} & A ,}
$$
$$
\xymatrix{ A \ot C \ar[d]_-{\lambda} \ar[r]^-{A \ot
\delta_C} & A\ot C \ot C \ar[r]^{\lambda \ot C} &
C\ot A \ot C \ar[d]^{C \ot  \lambda} & A\ot A
\ot C \ar[d]_{m_A \ot C} \ar[r]^{A \ot \lambda} & A
\ot C \ot A \ar[r]^{\lambda \ot A} & C \ot A
\ot A \ar[d]^{C \ot m_A}\\
C \ot A \ar[rr]_{\delta_C \ot A}&& C \ot C \ot A,
& A \ot C \ar[rr]_{\lambda} && C \ot A \,\,.}
$$

For any entwining  $(\uA, \uC,\lambda)$, the natural
transformation
$$\lambda'= \lambda\ot -:{\rm{T}}_{\uA}\circ {\rm{G}}^{\uC }=A \ot C \ot
- \to C \ot A \ot -=\rm{G}^{\uC }\circ {\rm{T}}_{\uA}$$
is a mixed distributive law from the monad
$\rm{T}_{\uA}$ to the comonad $\rm{G}^{\uC}$.
We write $\widetilde{\uC}$ for the ${_A \V}$-comonad
$\widetilde{\rm{G}^{\uC}}$, that is, for any $(V,h_V)\in
{_A\V}$,
$$\widetilde{\uC }(V,h_V)=(C\ot V, A\ot C \ot V
\xrightarrow{\lambda \ot V}C\ot A \ot V \xrightarrow{C\ot h_V} C \ot V),$$
and write ${_A^C\V }(\lambda)$ for the
category
$\V^{\rm{G}^{\uC} }_{\rm{T}_{\uA}}(\lambda')$. An
object of this category is a three-tuple $(V, \theta_V, h_V)$,
where $(V, \theta_V) \in {{^C\V}}$  and $(V, h_V) \in {_A\!\V}$,
with commuting diagram
\begin{equation}\label{mixed}
 \xymatrix{A \ot V \ar[r]^-{h_V} \ar[d]_{A \ot \theta_V}&
 V \ar[r]^-{\theta_V}& C \ot V \\
  A\ot C \ot V \ar[rr]_-{ \lambda \ot V}& & C \ot A \ot
  V\ar[u]_{C \ot h_V}.}
\end{equation}
 The assignment $((V,h_V),\theta_{(V,h_V)})\longmapsto
(V,\theta_{(V,h_V)},h_V)$
yields an isomorphism of categories
$$\Lambda : ({_A\V})^{\widetilde{\uC} } \to {_A^C\V }(\lambda).$$

\begin{thm}\label{rep-entw}{\bf Representable entwinings.} \em
For objects $A$, $C$ in a
monoidal category $\cV=(\V, \ot, \II)$, consider the functor $$\V(-\ot C, A): \V^{op}\to \rm{Set}$$
 taking an arbitrary object $V \in \V$ to the set $\V( V\ot C, A).$
Suppose there is an object $E\in \V$ that
represents the functor, i.e. there is a natural bijection $$\omega:
\V(-\ot C, A)\simeq \V(-,E).$$ Writing $\beta: E \ot C \to
A$ for the morphism $\omega^{-1}(\id_E)$, it follows that for any
object  $V \in \V$ and any morphism $f: V \ot C \to A$, there
exists a unique $\beta_f : V \to E$ making the diagram
$$
\xymatrix{ V \ot C \ar[rr]^{f} \ar[rd]_{\beta_f \ot C}&& A\\
   & E \ot C \ar[ru]_{\beta}  }$$
commute.
It is clear that $\omega^{-1}(\beta_f)=f.$

We call an entwining $(\uA,\uC,\lambda)$
\emph{representable} if the functor $\V(-\ot \!C, A):
\V^{op}\to \rm{Set}$ is representable.
\end{thm}

\begin{thm}\label{Ex.2}{\bf Examples.}\em
\begin{rlist}
    \item If the functor $-\ot C : \V \to \V$ has a right
    adjoint $[C,-]: \V \to \V$, then it follows from the
    bijection $\V(V \ot C, A) \simeq \V(V, [C,A])$ that the
    object $[C,A]$ represents the functor $\V(-\ot C, A)$. In
    particular, when $\cV$ is right closed, each entwining in
    $\cV$ is representable.
\item If $C$ is a right $\cV$-nuclear object, the functor
$\V(-\ot C, A): \V^{op}\to\rm{Set}$ is representable.
Indeed, to say that $C$ is right
$\cV$-nuclear is to say that the functor $C \ot - :\V \to \V$
has a right adjoint $\{C,-\}: \V \to \V$, the functor $-\ot
{^*C}=-\ot \{C,\II\}: \V \to \V$ also has a right adjoint
$[{^*C},-]: \V \to \V$, and the morphism $$\alpha'_V : V \ot
C \to [{^*C},V]$$ is a natural isomorphism. Considering then the
composition of bijections
$$\V(-\ot {^*C},?)\simeq \V(-, [{^*C},?]) \simeq \V(-, ? \ot
\!C),$$
one sees that the functor $-\ot C$ admits the functor
$-\ot {^*C}$ as a left adjoint. The same arguments as in
the proof of Theorem X.7.2 in \cite{M} then show  that there is a
natural bijection $$\V(-\ot \!C,A)\simeq \V(-, A \ot
{^*C} ).$$ Thus the object $A\ot {^*C}$ represents
the functor $\V(-\ot C, A): \V^{op}\to \rm{Set}$.
\end{rlist}
\end{thm}

\begin{proposition} \label{attached-alg}
 Let $(\uA,\uC,\lambda)$ be a representable entwining
in $\cV=(\V, \ot, \II)$
with representing object $E$ (see \ref{rep-entw}).
Let $e_E=\beta_\tau: \II \to E$ and $m_E=\beta_\varrho : E \ot
E \to E$, where $\tau$ is the composite $$\II \ot C\simeq C
\xrightarrow{\varepsilon_C}\II\xrightarrow{e_A}A,$$ while
$\varrho$ is the composite $$E\ot E \ot C
\xrightarrow{E\ot E \ot \delta_C} E \ot E\ot C
\ot C \xrightarrow{E\ot \beta \ot  C} E \ot A
\ot  C \xrightarrow{E \ot \lambda}E \ot C \ot A
\xrightarrow{\beta \ot A} A \ot  A \xrightarrow{m_A} A.$$
 \begin{rlist}
 \item The triple $(E, e_E,m_E)$ is a $\cV$-algebra.
 \item The morphism $i:=\beta_{A \ot \varepsilon_C}: A \to E$
 is a morphism of $\cV$-algebras.
 \end{rlist}
\end{proposition}

\begin{proof} (i) First observe commutativity of the diagrams
\begin{equation}\label{alg}
\xymatrix{E \ot E \ot C \ar[d]_{m_E \ot C}
\ar[rd]^{\varrho} & & \II \ot C=C \ar[d]_{e_E \ot C}
\ar[rd]^{e_A \cdot\varepsilon_C}&\\
E \ot C \ar[r]_-{\beta}& A, &  E \ot C \ar[r]_-{\beta}& A.}
\end{equation}

  To prove that $m_E$ is associative, i.e.
$m_E \cdot (m_E\ot E)=m_E \cdot (E \ot m_E)$,
it is  to show
$$\beta \cdot (m_E \ot C)
\cdot (m_E \ot E \ot C)=\beta \cdot (m_E \ot C)
\cdot(E \ot m_E \ot C).$$
 For this, consider the diagram (deleting the $\ot$-symbols)
$$\xymatrix{EEEC \ar[r]^-{EEE\delta_C}\ar[d]_{m_EEC}& EEECC \ar[r]^-{EE\beta
C}\ar[d]|{m_EECC}& EEAC \ar[r]^-{EE\lambda}\ar[d]|{m_ECA}&EECA
\ar[r]^-{EE\delta_C A}\ar[d]^{m_ECA}&EECCA \ar[r]^{E\beta CA}&EACA \ar[d]^{E \lambda A} \\
EEC \ar@{}[ru]|{(1)}\ar[r]_-{EE\delta_C}&EECC
\ar@{}[ru]|{(2)}\ar[r]_-{E\beta C}&EAC \ar@{}[ru]|{(3)}
\ar[r]_-{E\lambda}& ECA \ar[d]_{\beta A}&& ECAA \ar[d]^{\beta AA}\\
&&&AA \ar[d]_{m_A} \ar@{}[rruu]|{(4)}&&AAA \ar[d]^{Am_A}\ar[ll]_-{m_AA}\\
&&&A \ar@{}[rru]|{(5)}&&AA \ar[ll]^-{m_A}\,,}$$ in which
 the diagrams (1),(2) and (3) commute by naturality of
    composition,
  diagram (4) commutes by definition of $m_E$, and
  diagram (5) commutes by associativity of $m_A$.
It  follows
\begin{equation}\label{ass}
\beta \cdot (m_E  C) \cdot (m_E  E  C)=m_A
\cdot Am_A
 \cdot \beta AA \cdot E \lambda A \cdot  E\beta CA \cdot EE\delta_C A
 \cdot EE\lambda \cdot EE\beta C  \cdot EEE\delta_C.
\end{equation}

Now we have
$$\begin{array}{rcl}
\beta \cdot m_E C\cdot Em_EC&=&m_A \cdot \beta A \cdot E \lambda \cdot E\beta C \cdot EE\delta_C \cdot Em_EC\\[+1mm]
 \text{\scriptsize{nat. of composition}}&=
  &m_A \cdot \beta A \cdot E \lambda \cdot E\beta C \cdot
   Em_ECC \cdot EEE\delta_C  \\[+1mm]
\text{\scriptsize{by }(\ref{alg})}&=& m_A \cdot \beta A \cdot E \lambda \cdot Em_A C \cdot E\beta AC \cdot EE\lambda C \cdot EE \beta CC \cdot EEE\delta_C C \cdot EEE\delta_C\\[+1mm]
 \text{\scriptsize{coassociativity of } $\delta_C$} &=& m_A \cdot \beta A \cdot E \lambda
 \cdot Em_A C \cdot E\beta AC \cdot EE\lambda C \cdot EE \beta CC \cdot EEE C\delta_C  \cdot EEE\delta_C\\[+1mm]
 \text{\scriptsize{nat. of composition}}&=&m_A \cdot \beta A \cdot E \lambda
 \cdot Em_A C \cdot E\beta AC \cdot EE\lambda C \cdot EEA\delta_C \cdot EE\beta C  \cdot EEE\delta_C\\[+1mm]
 \text{$\lambda$ \scriptsize{is an entwining}}&=&m_A \cdot \beta A
 \cdot ECm_A\cdot E \lambda A \cdot  E A\lambda  \cdot  E\beta AC \cdot EE\lambda C \cdot EEA\delta_C \cdot EE\beta C  \cdot EEE\delta_C\\[+1mm]
 \text{\scriptsize{nat. of composition}}&=&m_A \cdot \beta A
 \cdot ECm_A\cdot E \lambda A \cdot  E\beta CA \cdot EEC\lambda \cdot EE\lambda C \cdot EEA\delta_C \cdot EE\beta C  \cdot EEE\delta_C\\[+1mm]
 \text{$\lambda$ \scriptsize{is an entwining}}&=&m_A \cdot \beta A
 \cdot ECm_A\cdot E \lambda A \cdot  E\beta CA \cdot EE\delta_C A
 \cdot EE\lambda \cdot EE\beta C  \cdot EEE\delta_C\\[+1mm]
 \text{\scriptsize{nat. of composition}}&=&m_A \cdot Am_A
 \cdot \beta AA \cdot E \lambda A \cdot  E\beta CA \cdot EE\delta_C A
 \cdot EE\lambda \cdot EE\beta C  \cdot EEE\delta_C\\[+1mm]
 \text{\scriptsize{by }(\ref{ass})}&=&\beta \cdot m_EC \cdot m_EEC.
\end{array}$$
It follows that $m_E \cdot (m_E \ot E)=m_E \cdot (E \ot m_E)$, 
and thus $m_E$ is associative.
To prove that $e_E$ is the unit for the multiplication $m_E$, it is
 to show
\begin{center}
$\beta \cdot (m_E \ot C)\cdot(e_E  \ot E \ot C)=\beta$ \quad and \quad
$\beta \cdot (m_E \ot C)\cdot(E  \ot e_E \ot C)=\beta.$
\end{center}
In the diagram
$$\xymatrix{
EC \ar@{}[rdd]|{(1)}\ar[dd]_{e_E EC}\ar[r]^-{E\delta_C}& ECC
\ar@{}[rdd]|{(2)} \ar[dd]|{e_EECC }\ar[r]^-{\beta C}& AC
\ar@{}[rdd]|{(3)} \ar[dd]|{e_EAC} \ar[r]^-{\lambda}& CA
\ar@{}[rdd]|{(4)} \ar[dd]|{e_ECA} \ar[r]^-{\varepsilon_C A}&
A\ar[dd]|{e_EA}\ar@{=}[ddr]&\\\\
EEC \ar@{}[rdd]|{(5)}\ar[r]_-{EE\delta_C}&EECC
\ar@{}[rdd]|{(6)}\ar[r]_-{E\beta C}& EAC \ar[r]_-{E
\lambda}& ECA \ar[r]_-{\beta A}& AA \ar[r]_-{m_A}& A\\\\
EC \ar[uu]^{Ee_EC} \ar[r]_-{E\delta_C}& ECC \ar[uu]|{Ee_E
CC}\ar[r]_-{E \varepsilon_C C}& EC \ar[uu]|{Ee_AC}\ar[uur]|{ECe_A}
\ar[r]_-{\beta}&A \ar@{}[uu]|{(7)}\ar[uur]|{Ae_A} \ar@{=}[r]&A\,,
\ar@{}[uu]|{(8)} \ar@{=}[uur]&}$$
  the diagrams (1), (2), (3), (5) and (7) are commutative by
    naturality of composition,
   the diagrams (4) and (6) are commutative by (\ref{alg}),
    diagram (8) and the top triangle are commutative since $e_A$ is the unit for
    the multiplication $m_A$, and
    the bottom triangle commutes since $\lambda$ is an entwining.
 Commutativity of the diagram implies
$$\begin{array}{rcl}
 \beta \cdot (m_E \ot C)\cdot(e_E  \ot E \ot C)
&=&(\varepsilon_C \ot A)\cdot \lambda \cdot(\beta \ot C)(E
\ot \delta_C) \\[+1mm]
\text{$\lambda$ \scriptsize{is an entwining}}&=&(A \ot
\varepsilon_C) \cdot(\beta \ot C)\cdot(E
\ot \delta_C) \\[+1mm]
\text{\scriptsize{nat. of composition}}&=&\beta \cdot (E
\ot C \ot\varepsilon_C)\cdot(E
\ot \delta_C) \\[+1mm]
\text{ \scriptsize{since}$(C \ot  \varepsilon_C)\cdot
\delta_C=\id$}&=&\beta \\[+1mm]
&=&\beta \cdot (m_E \ot C)\cdot(E  \ot e_E \ot C).
\end{array}$$
Thus $e_E$ is the unit for $m_E$. This completes the proof of (i).
\smallskip

(ii) We have to show $i\cdot e_A = i_E$ and
$m_E \cdot (i \ot i)=i \cdot m_A$.
For this consider the diagram
$$
\xymatrix{ C \ar[d]_{\varepsilon_C} \ar[r]^{e_A \ot C\quad}&
A \ot C \ar[r]^{i \ot C}\ar[d]_{A \ot \varepsilon_C} &
  E \ot C   \ar[ld]^{\beta} \\
\II\ar[r]_{e_A}  &A}$$
Since the square commutes by
naturality of composition, while the triangle commutes by
definition of $i$, the outer diagram is also commutative,
meaning  $\beta_{e_A \cdot \varepsilon_C}=i\cdot e_A.$
But $\beta_{e_A \cdot \varepsilon_C}=e_E.$ Thus $i\cdot e_A=e_E$.

Next, consider the diagram
$$
\xymatrix{& A \ot A \ot C \ot C \ar[rr]^-{A \ot A
\ot \varepsilon_C \ot C} \ar[d]^-{i \ot A \ot C\ot C}
& & A \ot A \ot C \ar[dd]^{i \ot A \ot C}
\ar[r]^{A \ot \lambda}&
A \ot C \ot A \ar[dd]_{i \ot C \ot A}\ar[ddr]^{A \ot \varepsilon_C \ot A}&\\
A \ot A \ot C \ar[ru]^{A \ot A \ot\delta_C}
\ar[d]_{i\ot i \ot C}& E \ot A \ot C\ot C
\ar[rrd]^{\;E \ot A \ot \varepsilon_C \ot C}\ar[d]_{E \ot i \ot C \ot C}&& &&\\
E \ot E \ot C \ar[r]_-{E \ot E \ot \delta_C}& E\ot E \ot C \ot C
\ar[rr]_-{E \ot \beta \ot C}
& & E\ot A \ot C \ar[r]_-{E \ot \lambda}& E \ot C \ot A
\ar[r]_-{\beta \ot A} & A \ot A }$$
in which the two triangles commute by definition of $i$, and the quadrangles
commute by naturality of composition. We have

$$\begin{array}{rcl}
\omega^{-1}(m_E \cdot (i \ot i))&=&m_A \cdot (\beta \ot A
)\cdot (E \ot \lambda) \cdot (E \ot \beta \ot C) \cdot
(E \ot E \ot \delta_C) \cdot (i\ot i \ot C)\\
&=&m_A \cdot (A \ot \varepsilon_C \ot  A) \cdot (A \ot
\lambda) \cdot (A \ot A\ot \varepsilon_C  \ot C)
\cdot (A \ot A \ot \delta_C)\\
\text{\scriptsize{since $(\varepsilon_C\ot C) \cdot \delta_C=\id$}}&=&m_A \cdot (A\ot \varepsilon_C \ot A)\cdot (A  \ot\lambda)  \\
\text{\scriptsize{by definition of $\lambda$}}
&=&m_A \cdot (A\ot A\ot \varepsilon_C) \\
 \text{\scriptsize{nat. of composition}} &=&
(A \ot \varepsilon_C)\cdot (m_A \ot C)   \\
\text{\scriptsize{since $A \ot\ve_C=\beta \cdot( i \ot C)$}}
         &=&\beta \cdot (i \ot C) \cdot (m_A \ot C) \\
&=&\beta \cdot ((i \cdot m_A)\ot C)\\
&=&\omega^{-1}(i \cdot m_A).
\end{array} $$
Thus $m_E \cdot (i \ot i)=i \cdot
m_A$. This completes the proof.
\end{proof}

\begin{proposition}\label{p-triangle} 
With the data given in Proposition \ref{attached-alg},
there is a functor
$$\Xi: {^C_A\V}(\lambda)\to
{_E\V}$$  with commutative diagram
\begin{equation}\label{triangle}
\xymatrix{^C_A\V(\lambda) \ar[rr]^{\Xi} \ar[dr]_{^C_AU}&& _E\V
\ar[dl]^{_EU}\\
&\V\,,&}\end{equation} where $^C_AU :{^C_A\V}(\lambda)\to \V$ is
the evident forgetful functor.
\end{proposition}
\begin{proof} For any $(V,\theta_V,h_V)\in {^C_A\V}(\lambda)$,
write $\iota_V : E \ot V \to V$ for the composite $$E \ot
V \xrightarrow{E \ot \theta_V}E \ot C \ot V
\xrightarrow{\beta \ot V}A \ot V \xrightarrow{h_V} V.$$ We
claim that $(V,\iota_V) \in {_E\V}$. Indeed, to show that $\iota_V
\cdot (e_E \ot V)=1,$ consider the diagram
$$
\xymatrix{ V \ar[d]_{e_E \ot V}\ar[r]^{\theta_V} & C \ot V
\ar[d]_{e_E \ot C \ot V}\ar[r]^-{\varepsilon_C
\ot V}& V \ar[d]_{e_{A} \ot V}  \ar@{=}[rd]&\\
E \ot V   \ar[r]_{E \ot\theta_V}& E \ot C \ot V
\ar[r]_{\beta \ot V}& A \ot V \ar[r]_{h_V}& V\,,}$$ in
which the left squares commutes by naturality of composition, the
right one commutes by (\ref{alg}), while the triangle
commutes since $(V,h_V)\in {_A\cV}$. It follows -- since
$(V,\theta_V) \in {^C\V}$, and hence $(\varepsilon_C \ot
V)\cdot \theta_V=\id$ -- that $\iota_V \cdot (e_E \ot
V)=(\varepsilon_C \ot V)\cdot \theta_V=\id.$

Next, consider the diagram
$$\xymatrix{EEV
\ar@{}[rrd]|{(1)}\ar[rr]^-{EE\theta_V}\ar[d]_-{EE\theta_V}&&
EECV \ar@{}[rrd]|{(2)}\ar[rr]^-{E\beta V}
\ar[d]_-{EEC\theta_V}&& EAV \ar@{}[rrd]|{(3)}
\ar[r]^-{Eh_V}\ar[d]_-{EA\theta_V}& EV \ar[rd]^-{E\theta_V}&\\
EECV \ar@{}[rrrrrdd]|{(4)}\ar[dd]_{m_ECV}\ar[rr]_-{EE\delta_C V}&&
EECCV \ar[rr]_-{E\beta CV}&& EACV
\ar[r]_-{E\lambda V}& ECAV \ar[d]_{\beta AV}\ar[r]^-{ECh_V}&ECV \ar[d]^{\beta V} \\
&&&&& AAV \ar@{}[ru]|{(5)}\ar[d]_{m_AV}\ar[r]_{Ah_V}&AV \ar[d]^{h_V}\\
ECV \ar[rrrrr]_{\beta V}&&&&& AV
\ar@{}[ru]|{(6)}\ar[r]_{h_V}&V\,,}$$ in which
  diagram (1) commutes since $(V,\theta_V) \in {^C\V}$,
 the diagrams (2) and (5) are commutative by naturality
    of compositions,
    diagram (3) commutes by (\ref{mixed}),
     diagram (4) commutes by definition of $m_E$, and
     diagram (6) commutes since $(V,h_V)\in {_A\V}$.
It follows that
$$\begin{array}{rcl}
\iota_V \cdot (E \ot \iota_V)&=&
 h_V \cdot (\beta \ot V ) \cdot (E \ot  \theta_V)\cdot (E
\ot h_V)\cdot (E \ot \beta \ot V)\cdot (E \ot  E\ot\theta_V)\\
&=&h_V \cdot (\beta \ot V ) \cdot (m_E  \ot  C \ot
V)\cdot (E \ot E \ot \theta_V) \\
\text{\scriptsize{nat. of composition}}
&=&h_V \cdot (\beta \ot V ) \cdot (E \ot \theta_V)\cdot(m_E  \ot V) \\
&=&\iota_V \cdot (m_E \ot V ).
\end{array}$$
 Thus, $(V,\iota_V)\in {_E\V}.$
Next, if $f: (V, \theta_V,h_V)\to (V',\theta_{V'},h_{V'})$ is
in ${^C_A\V}(\lambda)$, then the diagram
\begin{equation}\label{morph}
\xymatrix{ A \ot V \ar[d]_{A \ot f}\ar[r]^{h_V} & V
\ar[d]_{f}\ar[r]^-{\theta_V}& C \ot V \ar[d]^{C \ot f}\\
A \ot V'   \ar[r]_{h_{V'}}&  V' \ar[r]_{\theta_{V'}}& C
\ot V'}\end{equation} 
is commutative. In the diagram
$$
\xymatrix{ E \ot V \ar[d]_{E \ot f}\ar[r]^-{E \ot
\theta_V} & E \ot C \ot V
\ar[d]_{E \ot C \ot f}\ar[r]^-{\beta \ot V}& A \ot V \ar[d]_{A \ot f}\ar[r]^{h_V}& V \ar[d]^{f}\\
E \ot V'   \ar[r]_-{E \ot \theta_{V'}}&  E \ot C
\ot V' \ar[r]_{\beta \ot V'}& A \ot V'
\ar[r]_{h_{V'}}& V'\,,}$$ the middle square commutes by
naturality of composition, while the other squares commute by
(\ref{morph}). Thus, $f$ can be seen as morphism in ${_E\V}$ from $(V,\iota_V)$ to
$(V',\iota_{V'})$. It follows that the assignment
$$(V, \theta_V,h_V)\longmapsto (V, \iota_{V}=h_V \cdot (\beta
\ot V)\cdot (E \ot \theta_V))$$ yields a functor
$\Xi: {^C_A\V}(\lambda)\to {_E\V}.$
  It is clear that $\Xi$ makes the diagram (\ref{triangle}) commute.
\end{proof}

In order to proceed, we need the following result.

\begin{lemma}\label{l-square}
Let $\bT=(T,e_T,m_T)$ and $\textbf{H}=(H,e_H,m_H)$ be
monads on a category $\A$, $i:T \to H$ a monad morphism and
$i_*:\A_H \to \A_T$  the functor that takes an $H$-algebra
$(a,h_A)$ to the $T$-algebra $(a, h_a \cdot i_a).$ Suppose that
$\bT^\di=(T^\di,e^\di ,m^\di )$ (resp.
$\textbf{H}^\di=(H^\di,e_{H^\di},m_{H^\di})$) is a comonad that is
right adjoint to $\bT$ (resp. $\textbf{H}).$ Write
$\overline{i}: H^\di \to T^\di$ for the mate of $i$. Then

\begin{zlist}
\item $\overline{i}$ is a morphism of comonads.
\item We have commutativity of the diagram
\begin{equation}\label{square}
\xymatrix{ \A_H \ar[d]_{K_{H,H^\di}}\ar[r]^-{i_*} &  \A _T\ar[d]^{K_{T,T^\di}}\\
{\A}^{H^\di} \ar[r]_{(\overline{i})_*}& {\A}
^{T^\di}.}\end{equation}

\item If $\A$ admits
both equalisers and coequalisers, then 
 \begin{rlist}
\item the functors $i_*$ and
$(\overline{i})_*$ admit both right and left adjoints; 
\item 
$i_*$ is monadic and $(\overline{i})_*$ is comonadic.
\end{rlist}
\end{zlist}
\end{lemma}

\begin{proof} 
(1) This follows from the properties of mates  (see, for example, \cite{MW}).

(2) An easy calculation shows that for any $(a,h_a)\in
\A_H$, one has
$$\begin{array}{ll}
(\overline{i})_* \circ K_{H,H^\di}(a,h_a)=(a,
(\overline{i})_a \cdot H^\di(h_a) \cdot \sigma_a) & \text{and } \\[+1mm]
K_{T,T^\di}\circ i_*(a,h_a)=(a, T^\di(h_a) \cdot T^\di(i_a) \cdot \tau_a),
\end{array}$$
where $\sigma : \id \to H^\di H$ is the unit of the adjunction $H
\dashv H^\di,$ while $\tau : \id \to T^\di T$ is the unit of the
adjunction $T \dashv T^\di.$ Considering the diagram
$$\xymatrix{a \ar[r]^-{\sigma_a} \ar[d]_{\tau_a}& H^\di H (a)
\ar[d]^{\overline{i}_{H(a)}}\ar[r]^-{H^\di(h_a)} & H^\di(a) \ar[d]^{(\overline{i})_a}\\
T^\di T(a) \ar[r]_-{T^\di(i_a)}& T^\di H(a) \ar[r]_-{T^\di (h_a)}&
T^\di(a)\,,}$$ in which the right square commutes by naturality of
$\overline{i}$, while the left one commutes by Theorem IV.7.2 of
\cite{M}, since $\overline{i}$ is the mate of $i$. Thus $
(\overline{i})_a \cdot H^\di(h_a) \cdot \sigma_a=T^\di(h_a) \cdot
T^\di(i_a) \cdot \tau_a$, and hence $(\overline{i})_* \circ
K_{H,H^\di}=K_{T,T^\di}\circ i_*$.

(3)(i) If the category $\A$ admits both equalisers and coequalisers, then
$\A_H$ admits equalisers, while $A^{H^\di}$ admits coequalisers.
Since $K_{H,H^\di}$ is an isomorphism of categories, it follows
that the category $\A_H$ as well as the category $A^{H^\di}$ admit
both equalisers and coequalisers. Then, according to
\ref{right-ad} and its dual, the functor $(\overline{i})_*$ admits
a right adjoint $(\overline{i})^*$, while the functor $i_*$ admits
a left adjoint $i^*.$ Then clearly the composite
$i_!=(K_{H,H^\di})^{-1}\cdot (\overline{i})^* \cdot K_{T,T^\di}$
is right adjoint to $i_*$, while the composite
$(\overline{i})_!=K_{H,H^\di}\cdot i^* \cdot (K_{T,T^\di})^{-1}$
is left adjoint to $(\overline{i})_*$.

(3)(ii) Since the functors $i_*$ and
$(\overline{i})_*$ are clearly conservative, 
the assertion follows by a simple application of Beck's
monadicity theorem  (see, \cite{M}) and its dual.
\end{proof}

\begin{thm}\label{left-right-adj}{\bf Left and right adjoints to the functor $i_*$.} \em In the setting considered in Proposition \ref{attached-alg},
suppose now that $\V$ admits both equalisers and coequalisers and
that there are adjunctions $A\ot - \dashv \{A,-\}: \V \to \V$
and $E\ot - \dashv \{E,-\}: \V \to \V.$ Then, according to  
 Lemma \ref{l-square}, the functor $i_* :{_E\!\V} \to {_A\!\V}$ has both left
and right adjoints $i^*$ and $i_!$. It follows from \ref{right-ad}
and its dual that for any $(V, h_V)\in {_A\!\V}$, $i^*(V, h_V)$ is
the coequaliser
$$\xymatrix{E \ot A  \ot V
\ar@<-0.5ex>[rr]_-{E \ot h_V} \ar@<0.5ex>[rr]^-{h^r_E \ot V}
&& E \ot V \ar[r]^{q_{(V, h_V)}} & i^*(V, h_V)\,,}$$ while
$i_!(V, h_V)$ is the equaliser
\begin{equation}\label{equaliser}
\xymatrix{i_!(V, h_V) \ar[r]^-{e_{(V, h_V)}}  &\{E,V\}
\ar@<-0.5ex>[rr]_-{k} \ar@<0.5ex>[rr]^-{\{h^l_E,V\}} & & \{A \ot E,V\}\,,}
\end{equation}
where $h^r_E=m_E\cdot ( E\ot i)$, $h^l_E=m_E\cdot ( i\ot E)$ and $k$ is the transpose of the composition
 $$A\ot E \ot \{E,V \}\xrightarrow{A\ot \overline{e}^E_V }A
\ot V\xrightarrow{h_V}V.$$
We write $E \ot_A -$ for the functor $i^*$ (as well as for the ${_A\!\V}$-monad generated by the adjunction $i^* \dashv i_*$) and
write $\{E,-\}_A$ for the functor $i_!$ (as well as for the
${_A\!\V}$-comonad generated by the adjunction $i_* \dashv i_!$).

According to Lemma \ref{l-square}, the comparison functor
$K_{i^*}:{_E\V} \to({_A\V})_{E\ot_A -}$ is
an equivalence of categories. It is easy to check that for any
$(V,h_V)\in {_E\V},$
$K_{i^*}(V,h_V)=((V,\nu_V),\kappa_{(V,h_V)})$, where $\nu_V=h_V
\cdot (i\ot V)$, while $\kappa_{(V,h_V)}:E \ot_A V \to V$
is the unique morphism with commutative diagram
\begin{equation}\label{triangle1}
\xymatrix{E \ot V \ar[rr]^{q_{(V,\nu_V)}} \ar[dr]_{h_V}&&
E\ot_A V \ar[dl]^{\kappa_{(V,h_V)}}\\
&V&.}\end{equation}
Such a unique morphism exists because
the morphism $h_V : E \ot V\to V$ coequalises the pair of
morphisms $(h_E\ot V,E\ot \nu_V ).$
\end{thm}
\begin{thm}\label{pairing}{\bf Pairing induced by the adjunction $E \ot_A -\dashv \{E,-\}_A.$} \em
We refer to the setting considered in  \ref{left-right-adj}.
Writing $\mathfrak{F}$ for the composition
$$({_A\V})^{\widetilde{\uC}}\xrightarrow{\Lambda }{_A^C\V }(\lambda)
\xrightarrow{\Xi}{_E\V} \xrightarrow{K_{i^*}}({_A\V})_{E\ot_A
-}\xrightarrow{K_{E\ot_A -,\,\{E,-\}_A}}(_A\V)^{\{E,-\}_A},$$
one easily sees that $\mathfrak{F}$ makes the diagram
$$\xymatrix{({_A\V})^{\widetilde{\uC}} \ar[rr]^{\mathfrak{F}} \ar[dr]_{U^{\widetilde{\uC}}}&&
({_A\V})^{\{E,-\}_A} \ar[dl]^{_{\{E,-\}_A}U}\\
&{_A\V}&}$$ commute, were $U^{\widetilde{\uC}}$ is the
evident forgetful functor. Then, according to \ref{commorphism},
there is a unique morphism of comonads
$\alpha:\widetilde{\uC} \to \{E,-\}_A$ such that
$\alpha_*=\mathfrak{F}$.

Since the triple $(E \ot_A -, \{E,-\}_A, \widehat{e}^E_{-}) $, where  $\widehat{e}^E_{-}$
is the counit of the adjunction $E \ot_A - \dashv \{E,-\}_A$, is a pairing (see \ref{pair-func}),
it follows from Proposition
\ref{Example 2} that the triple 
\begin{equation}\label{p-lambda}
\pp(\lambda)=(E\ot_A-,\{E,-\}_A,\sigma:=\widehat{e}^E_{-}\cdot(E \ot_A\alpha))
\end{equation}
is a pairing on the category ${_A\!\V}$.   

A direct inspection shows that for any
$((V,\nu_V),\theta_{(V,\nu_V)})
\in({_A\V})^{\widetilde{{\uC}}}$, $\Xi\Lambda
((V,\nu_V),\theta_{(V,\nu_V)})=(V,\xi),$ where $\xi$ is the
composite
$$E \ot V \xrightarrow{E \ot \theta_{(V,\nu_V)}} E \ot C \ot
V \xrightarrow{ \beta \ot V} A \ot V \xrightarrow{\nu_V}
V.$$

Then $$K_{i^*}\Xi\Lambda
((V,\nu_V),\theta_{(V,\nu_V)})=K_{i^*}(V,\xi)=((V,\nu_V),\kappa_{(V,\xi)}:
E \ot_A V \to V),$$ and thus
$$\mathfrak{F}((V,\nu_V),\theta_{(V,\nu_V)})=K_{E\ot_A -,\,\{E,-\}_A}((V,\nu_V),
\kappa_{(V,\xi)})=((V,\nu_V),\overline{\theta}_{(V,\nu_V)}),$$
where $\overline{\theta}_{(V,\nu_V)}$ is the composite
$$V \xrightarrow{\widehat{\eta}^E_V}\{E, E \ot_A V\}_A
\xrightarrow{\{E,\kappa_{(V,\xi)}\}_A}\{E,V\}_A.$$ 
Here $\widehat{\eta}^E_{-}: \id \to \{E,E\ot_A -\}_A$ is the unit of
the adjunction $E \ot_A - \dashv \{E,-\}_A.$

Since for any object $(V,\nu_V)\in \V_A, $ the pair
$$(\widetilde{C}(V,\nu_V), (\delta_{\widetilde{C}})_{(V,\nu_V)})=
((C \ot V, h), \delta_C \ot V),$$ where $h$ is the
composite $(C \ot \nu_V )\cdot (\lambda \ot V):A \ot C
\ot V \to C \ot V$, is an object of the category
$(_A\V)^{\widetilde{C}}$, one has
$$\mathfrak{F}((C \ot
V, h), \delta_C \ot V)=((C \ot V, (C\ot \nu_V )\cdot
(\lambda \ot V),\overline{\theta}_{(C \ot V,h)}),$$ where
$\overline{\theta}_{(C \ot V,h)}: C \ot V \to \{E, E
\ot_A V\}_A$ is the composite
$$C \ot V \xrightarrow{\widehat{\eta}^E_{C \ot V}} \{E, E\ot_A (C \ot V)\}_A
\xrightarrow{\{E,\kappa_{(C \ot V,\xi)}\}}\{E, C \ot
V\}_A.$$ Here $\xi$ is the composite
$$E \ot C \ot V\xrightarrow{E \ot \delta_C \ot V}E \ot C  \ot C  \ot
V\xrightarrow{\beta \ot C  \ot V} A \ot C \ot V
\xrightarrow{\lambda \ot V} C \ot A  \ot
V\xrightarrow{C\ot \nu_V } C \ot V.$$

Now, according to \ref{commorphism}, the $(V,\nu_V)$-component
$\alpha_{(V,\nu_V)}$ of the comonad morphism $\alpha:
\widetilde{C} \to \{E,-\}_A$ is the composite
$$C \ot V \xrightarrow{\widehat{\eta}^E_{C \ot V}} \{E, E\ot_A (C \ot V)\}_A
\xrightarrow{\{E,\kappa_{(C \ot V,\xi)}\}}\{E, C \ot V\}_A
\xrightarrow{\{E, \varepsilon_C \ot V\}_A} \{E,V\}_A,$$ i.e.
$\alpha_{(V,\nu_V)}$ is the transpose of the composite
$(\varepsilon_C \ot V)\cdot \kappa_{(C \ot V,\xi)}$. By (\ref{triangle1}),
the diagram
$$\xymatrix{E \ot C  \ot V \ar[d]_{q_{(C \ot CV,\xi)}}\ar[rd]^{\xi}&&\\
**[l]E \ot_A (C \ot V) \ar[r]_-{\kappa_{(C \ot
V,\xi)}}& C \ot V \ar[r]_-{\varepsilon_C \ot V}&V}$$
commutes. In the diagram
$$\xymatrix{ **[l]E \ot C \ot V \ar[r]^-{E \ot \delta_c \ot
V}\ar@{=}[dr]& E \ot C \ot C \ot V \ar[r]^-{\beta
\ot C \ot V}\ar[d]^{E \ot \varepsilon_C \ot C
\ot V}& A \ot C \ot V \ar[dr]_{A \ot\varepsilon_C \ot C}
\ar[r]^-{ \lambda \ot V}& C \ot A\ot V
\ar[r]^-{C\ot \nu_V }\ar[d]^{ \varepsilon_C \ot A \ot V}
& C \ot V \ar[d]^{\varepsilon_C \ot V}\\
& E \ot  C \ot V \ar[rr]_-{\beta \ot V}&&  A \ot V
\ar[r]_{\nu_V}&V\,,}$$
the left triangle commutes since $\uC$ is a coalgebra in $\cV$,
 the middle triangle commutes since the triple
    $(\uA,\uC,\lambda)$  is an entwining, and
 the trapeze and the rectangle commute by naturality of composition,
 hence $(\varepsilon_C \ot V)\cdot \xi=\nu_V \cdot (\beta \ot V).$  
Thus, $\kappa_{(C \ot V, \xi)}$ is the unique morphism
that makes the diagram
\begin{equation}\label{action-coaction}
\xymatrix{E \ot C \ot V \ar[r]^-{\beta \ot V}
\ar[d]_{q_{(C \ot V, h)}}& A
\ot V \ar[d]^{\nu_V}\\
**[l]E \ot_A (C \ot V) \ar[r]_-{(\varepsilon_C \ot
V)\cdot\kappa_{(C \ot V, h)}}& V}
\end{equation} commute. (Recall that here $h=(C \ot \nu_V
)\cdot ( \lambda\ot V): A\ot C \ot V \to C \ot
V$.)

\begin{proposition} For any $(V,\nu_V)\in {_A\V}$, consider the
morphism $\alpha'_{(V,\nu_V)}: C \ot V  \to \{E,V\}$ that is
the transpose of the composition $E \ot C  \ot V
\xrightarrow{\beta \ot V} A \ot V \xrightarrow{\nu_V} V .$
Then the diagram
$$
\xymatrix{\{E,V\}_A \ar[rr]^-{e_{(V,\nu_V)}}&& \{E,V\}\\
& \ar[lu]^{\alpha_{(V,\nu_V)}} C \ot V
\ar[ru]_{\alpha'_{(V,\nu_V)}}&}$$ is commutative.
\end{proposition}
\begin{proof} We show first that
\begin{equation}\label{eq.5.8}
\{h_E,V\}\cdot \alpha'_{(V,\nu_V)}=k \cdot \alpha'_{(V,\nu_V)}
\end{equation} (see the equaliser diagram (\ref{equaliser})). Since
the transpose of the morphism $k$ is
the composite $$A \ot E \ot \{E,V\} \xrightarrow{A \ot
\overline{e}^E_V}A \ot V \xrightarrow{h_V} V,$$ while the
transpose of the morphism $\{h^l_E,V\}$ is
the composite
$$A \ot E \ot  \{E,V\} \xrightarrow{A \ot E \ot  \{h^l_E,V\}}A \ot E \ot  \{A \ot E,V\} \xrightarrow{\overline{e}^{A \ot E}
_V}V$$ and this is easily seen to be the composite
$$A \ot E \ot \{E,V\} \xrightarrow{h^l_E \ot \{E,V\}}E \ot
\{E,V\}\xrightarrow{\overline{e}^E_V}V.$$ Thus, since $$(h^l_E \ot
\{E,V\})\cdot (A\ot E \ot \alpha'_{(V,\nu_V)})=(E \ot
\alpha'_{(V,\nu_V)})\cdot (h^l_E \ot  C\ot  V)$$ by naturality of composition, it is
enough to show that the diagram
$$
\xymatrix{**[l] A \ot E \ot  C \ot V \ar[r]^-{A\ot
E \ot \alpha'_{(V,\nu_V)}} \ar[d]_{h^l_E \ot  C\ot  V }&
A \ot E \ot \{E,V\} \ar[r]^-{A \ot \overline{e}^E_V
}\ar@{-->}[d]_{h^l_E \ot \{E,V\}}&**[r]A\ot  V \ar[d]^{\nu_V}\\
**[l]E \ot  C\ot  V \ar[r]_-{E \ot
\alpha'_{(V,\nu_V)}}& E \ot \{E,V\}
\ar[r]_-{\overline{e}^E_V}& V}$$ is commutative. Since
$\overline{e}^E_V \cdot (E \ot \alpha'_{(V,\nu_V)})=
\nu_V \cdot (\beta\ot V)$,
the diagram can be rewritten as
\begin{equation}\label{sq}
\xymatrix{**[l]A \ot E \ot  C \ot V \ar[d]_{h^l_E
\ot  C\ot  V } \ar[r]^-{A \ot (\nu_V \cdot (\beta
\ot V))}
& A \ot V \ar[d]^{\nu_V}\\
**[l]E \ot C \ot V \ar[r]_{\nu_V \cdot (\beta \ot V)}&V\,.}
\end{equation}

Consider now the diagram
$$
\xymatrix{ A \ot E \ot  C \ar[r]^-{A \ot  E \ot
\delta_C} \ar[d]_{i \ot E \ot C}& A \ot E\ot
C\ot C \ar[r]^-{A \ot \beta \ot C} \ar[d]_{i \ot E
\ot C\ot C} & A\ot A\ot C \ar[r]^-{A \ot
\lambda} \ar[d]_{i \ot A\ot C}&A \ot C \ar[d]_{i
\ot C\ot A}\ot A \ar[rd]^{A \ot
\varepsilon_C \ot A}&\\
E\ot E\ot C \ar[r]_-{E\ot E\ot \delta_C}&E\ot
E\ot C\ot C \ar[r]_-{E \ot \beta \ot C}&E\ot
A\ot C \ar[r]_-{E \ot  \lambda}& E\ot C\ot A
\ar[r]_-{\beta \ot A}&A\ot A\,,}$$ in which the three
rectangles commute by naturality of composition, while the triangle
commutes since $\omega^{-1}(A\ot \varepsilon_C)=\beta \cdot (i
\ot C)$. Recalling now that $h^l_E=m_E \cdot (i \ot E),$ we
have
$$\begin{array}{rcl}
\beta \cdot (h^l_E \ot C)&=&\beta \cdot (m_E \ot C) \cdot (i\ot E \ot C)\\
&=&m_A \cdot (\beta \ot A)\cdot (E \ot \lambda)\cdot (E
\ot \beta \ot C) \cdot (E\ot E\ot \delta_C)\cdot(i\ot E \ot C)\\
&=& m_A \cdot (A \ot \varepsilon_C \ot A )\cdot (A \ot
\lambda)\cdot (A \ot \beta \ot C)\cdot (A\ot E\ot\delta_C)\\
 \text{\scriptsize{since $\lambda$ is an entwining}} &=&m_A \cdot (A \ot A \ot
 \varepsilon_C )\cdot (A \ot\beta \ot C)\cdot (A\ot E\ot \delta_C)  \\
 \text{\scriptsize{nat. of composition}} &=&m_A \cdot
(A \ot \beta)\cdot (A \ot E \ot C \ot\ve_C )\cdot (A\ot E\ot \delta_C) \\
 \text{\scriptsize{since $(C \ot \ve_C) \cdot\delta_C=\id_C$}} &=&m_A \cdot (A \ot \beta) .
\end{array}$$
 Therefore,  $m_A \cdot (A \ot \beta)=\beta \cdot (h^l_E \ot C)$,  and hence
$$(m_A \ot V)\cdot (A \ot \beta \ot
V)=(\beta \ot V) \cdot (h^l_E \ot C \ot V).$$
Using now
that $\nu_V \cdot (A \ot\nu_V)=\nu_V \cdot (m_A \ot V)$,
since $(V,\nu_V) \in {_A\!\V}$, one has
$$\nu_V \cdot ( A \ot \nu_V)\cdot (A\ot \beta \ot
V)=\nu_V \cdot (\beta \ot V) \cdot (h^l_E \ot C \ot
V).$$ Thus the diagram (\ref{sq}) commutes. It follows that
$\{h^l_E,V\}\cdot \alpha'_{(V,\nu_V)}=k \cdot
\alpha'_{(V,\nu_V)}$, and since the diagram (\ref{equaliser}) is
an equaliser, there exists a unique morphism $\gamma_{(V,\nu_V)}:
C \ot V \to \{E,V\}_A$ that makes the diagram
$$\xymatrix{\{E,V\}_A \ar[rr]^-{e_{(V,\nu_V)}}&& \{E,V\}\\
& \ar[lu]^{\gamma_{(V,\nu_V)}} C\ot V\ar[ru]_{\alpha'_{(V,\nu_V)}}&}$$
commute. We claim that
$\gamma_{(V,\nu_V)}=\alpha_{(V,\nu_V)}.$ To see this, consider the diagram
$$
\xymatrix{**[l] E \ot C \ot  V \ar[r]^-{E\ot
\gamma_{(V,\nu_V)}} \ar[d]_{q_{(C \ot V, h)}}& E\ot
\{E,V\}_A \ar[r]^-{E\ot
e_{(V,\nu_V)}}\ar[d]_{{q_{\{E,V\}_{\!A}}}}&
**[r]E\ot \{E,V\} \ar[d]^{\overline{e}^E_V}\\
**[l]E\ot_A (C \ot  V) \ar[r]_-{E \ot_A
\gamma_{(V,\nu_V)}}& E \ot_A \{E,V\}_A
\ar[r]_-{\widehat{e}^E_V}& V\,,}$$ where $\widehat{e}^E_-$ is the
counit of the adjunction $-\ot_A E \dashv [E,-]_A$. In this
diagram, the left rectangle commutes by naturality of $q$ (recall
that $\gamma_{(V,\nu_V)}:C \ot V \to \{E,V\}_A$ is a morphism
in ${_A\V}$), while the right one commutes by definition of
$\widehat{e}^E_-$. Since
$$\overline{e}^E_V \cdot (E\ot
e_{(V,\nu_V)}) \cdot (E\ot
\gamma_{(V,\nu_V)})=\overline{e}^E_V \cdot (E\ot
(e_{(V,\nu_V)} \cdot \gamma_{(V,\nu_V)}))=\overline{e}^E_V \cdot
(E\ot \alpha'_{(V,\nu_V)}),$$ and since $\overline{e}^E_V
\cdot (E\ot \alpha'_{(V,\nu_V)})=\nu_V \cdot
(\beta\ot V)$, it follows that the diagram
$$
\xymatrix{E \ot C \ot V \ar[r]^-{ \beta\ot V}
\ar[d]_{q_{(C \ot V, h)}}& A
\ot V \ar[d]^{\nu_V}\\
**[l]E \ot_A (V \ot C) \ar[r]_-{\widehat{e}^E_V \cdot (E
\ot_A \gamma_{(V,\nu_V)})}& V}
$$ commutes. Comparing this diagram with (\ref{action-coaction}),
one sees that $$\widehat{e}^E_V \cdot (E \ot_A
\gamma_{(V,\nu_V)})=(\varepsilon_C \ot V)\cdot\kappa_{(C
\ot V, h)}.$$ Thus $\gamma_{(V,\nu_V)} : C \ot V \to
\{E,V\}_A$ is the transpose of the morphism $(\varepsilon_C
\ot V)\cdot\kappa_{(C \ot V, h)}$. Thus
$\gamma_{(V,\nu_V)}$ is just $\alpha_{(V,\nu_V)}$. This completes
the proof.
\end{proof}
\end{thm}

When the pairing $\pp(\lambda)$ (\ref{p-lambda}) is rational, we write
$\R^{\pp(\lambda)}(E)$ for the full subcategory of the category
${_E\V}$ generated by those objects whose images under the functor
$K_{i^*}$ lie in the category $\R^{\pp(\lambda)}(E \ot_A-)$.

The following result extends \cite[Theorem 3.10]{Abu},
\cite[Proposition 2.1]{EG}, and \cite[Theorem 2.6]{EGL} 
from module categories to monoidal categories. 

\begin{theorem}\label{th-entw} Let $\cV=(\V,\ot,\II)$ be a monoidal category
with $\V$ admitting both equalisers and coequalisers, and
$(\uA,\uC,\lambda)$ a representable entwining with
representable object $E$. Suppose that
\begin{zlist}
    \item the functors $A \ot -,\,E \ot - : \V \to \V$
    have right adjoints $\{A,-\}$ and $\{E,-\}$,
 \item for any $(V,\nu_V)\in {_A\V}$, the transpose $\alpha'_{(V,\nu_V)}: C \ot V
 \to \{E,V\}$ of the composite \\
    $E \ot C \ot V \xrightarrow{\beta \ot V}A
    \ot  V \xrightarrow{\nu_V} V$ is a monomorphism, and
 \item \begin{rlist}
   \item the functor $C \ot - :\V \to \V$ preserves equalisers, or
   \item the category $\V$ admits pushouts and the functor $C \ot - :\V \to \V$
        preserves regular monomorphisms and has a right adjoint, or
   \item the category $\V$ admits pushouts, every monomorphism in $\V$ is regular
        and the functor $C \ot - :\V \to \V$ has a right adjoint.
 \end{rlist}
\end{zlist} 
Then the pairing $\pp(\lambda)$ is rational and
there is an equivalence of categories $${^C_A\V(\lambda)} \simeq
\R^{\pp(\lambda)}(E).$$
\end{theorem}
\begin{proof} Since $e_{(V,\nu_V)}:\{E,V\}_A \to\{E,V\}$ is an
equaliser for all ${(V,\nu_V)}\in {_A\V}$, $\alpha_{(V,\nu_V)}$ is
a monomorphism if and only if $\alpha'_{(V,\nu_V)}$ is so. Thus,
the pairing $\pp(\lambda)$ is rational if and only if for any
$(V,\nu_V)\in {_A\V}$, the morphism $\alpha'_{(V,\nu_V)}: C\ot
V \to \{E,V\}$ is a monomorphism. Thus, condition (2) implies
that the pairing $\pp(\lambda)$ is rational.

Next, since the forgetful functor ${_AU}:{_A\V} \to\V$ preserves
and creates equalisers, the functor $\widetilde{C}:{_A\V} \to
{_A\V}$ preserves equalisers if and only if the composite
${_AU}\widetilde{C}:{_A\V} \to\V$ does so. But for any
$(V,\nu_V)\in {_A\V}$, ${_AU}\widetilde{C}(V,\nu_V)=C \ot V$.
It follows that if the functor $C \ot -:\V \to\V$ preserves
equalisers, then the functor $\widetilde{C}:{_A\V} \to {_A\V}$
does so. Since each of the conditions in (3) implies that the functor $C
\ot - :\V \to \V$ preserves equalisers (see the proof of
Proposition \ref{P.3.4}) and since the functor $E\ot_A- : _A\!\!\V \to _E\!\!\V$
has a right adjoint $\{E,-\}_A : _A\!\!\V \to _E\!\!\V$, one can apply Theorem \ref{T.2.15} to
get the desired result.
\end{proof}

Since for any $\cV$-coalgebra $\bC=(C,\varepsilon_C,\delta_C)$,
the identity morphism $\id_C : C \ot \II=C \to C= \II \ot C$
is an entwining from the trivial $\cV$-algebra
$\mathcal{I}=(\II,\id_{\II}, \id_{\II})$ to the $\cV$-coalgebra $\bC$,
it follows from Example \ref{Ex.2}(1) that this entwining is
representable with representable object $C^*=[C,\II]$. Applying
Proposition \ref{attached-alg} gives:
\begin{thm}\label{bicl-coalg}{\bf Coalgebras in monoidal closed categories.}
Assume the monoidal category $\cV$ to be closed and consider any
$\cV$-coalgebra  $\bC=(C,\varepsilon_C,\delta_C)$. Then the triple
$$\bC^*=(C^*=[C,\II],e_{C^*}, m_{C^*})$$ is a $\cV$-algebra, where
$e_{C^*}=\pi(\varepsilon_C)$, while $m_{C^*}$ is the morphism $C^*
\ot C^* \to C^*$ that corresponds to the composite
$$C^*\ot C^* \ot C
\xrightarrow{C^*\ot C^* \ot \delta_C} C^*\ot C^*
\ot C \ot C \xrightarrow{C^*\ot  e^C_\II \ot
C}C^*\ot C \xrightarrow{e^C_\II}\II$$ under the bijection (see
(\ref{E.24}))
$$\pi=\pi_{C^*\ot C^*, C,\II}: \V(C^*\ot C^* \ot C, \II)\simeq \V(C^*\ot C^*,C^*).$$
\end{thm}

\begin{proposition}\label{coal-pairing} In the situation of \ref{bicl-coalg},
the triple
$$\pp(\bC)=(\bC^*,\bC,\,t=e_\II^C: C^*\ot C\to \II)$$ is a left pairing in $\cV$.
\end{proposition}
\begin{proof} We just note that the equalities  $\pi^{-1}(e_{C^*})=\varepsilon_C$  and
$$\pi^{-1}(m_{C^*})=e^C_\II \cdot (C^*\ot e^C_\II \ot
C)\cdot (C^*\ot C^* \ot \delta_C)$$ imply commutativity of
the diagrams
$$
\xymatrix{  C^*  \ot C^* \ot C \ar[rr]^-{ C^* \ot C^*
\ot \delta_{C}} \ar[d]_{ m_{C^*}\ot C}&&C^* \ot C^*
\ot C \ot C \ar[rr]^-{C^* \ot e^C_\II   \ot C}&&
C^* \ot C \ar[d]_{e^C_\II } & \ar[l]_-{e_{C^*}\ot   C} C
\ar[ld]^{\varepsilon_{C}} \\
C^* \ot C \ar[rrrr]_-{e^C_\II }&&&&\II.}
$$
\end{proof}
Applying now either Proposition \ref{P.3.7} or Proposition
\ref{coal-pairing} yields

\begin{theorem}\label{P.3.11}  Let
Let $\cV=(\V, \ot, \II, [-, -])$ be a monoidal closed
category, $\oC=(C,\varepsilon_C,\delta_C)$  a $\cV$-coalgebra with
$C$ $\cV$-prenuclear, and assume $\V$ to admit equalisers. If
either
\begin{rlist}
    \item the functor $C \ot -: \V \to \V$
          preserves equalisers, or
    \item $C \ot -: \V \to \V$ admits pushouts and the functor
          $C \ot -: \V \to \V$ preserves regular monomorphisms, or
    \item $\V$ admits pushouts and every monomorphism in $\V$ is regular,
\end{rlist}
then $\R^{\pp(\oC)}(\oC^*)$ is a full coreflective subcategory of
${_{\oC^*}\V}$ and the functor $\Phi^{\pp(\oC)}: {^{\oC}\V} \to
{_{\oC^*}\!\V}$ corestricts to an equivalence
$R^{\pp(\oC)}:{^{\oC}\V} \to  \R(\oC^*)$.
\end{theorem}

A special case of the situation described in Theorem \ref{P.3.11}
is given by a finite dimensional $k$-coalgebra $C$ over a field $k$
and $\V$ the category of $k$-vector spaces.

\smallskip

{\bf Acknowledgements.} The work on this paper was started during
a visit of the first author at the Department of Mathematics at
the Heinrich Heine University of D\"usseldorf supported by the
German Research Foundation (DFG) and continued with support by
Volkswagen Foundation  (Ref.: I/84 328). The authors express their
thanks to all these institutions.

\bigskip

\noindent
{\bf Addresses:} \\[+1mm]
{Razmadze Mathematical Institute, 1, M. Aleksidze st., Tbilisi
0193,  } {\small and} \\
 {Tbilisi Centre for Mathematical Sciences,
Chavchavadze Ave. 75, 3/35, Tbilisi 0168}, \\
 Republic of Georgia,
    {\small bachi@rmi.acnet.ge}\\[+1mm]
{Department of Mathematics of HHU, 40225 D\"usseldorf, Germany},
  {\small wisbauer@math.uni-duesseldorf.de}

\end{document}